\documentclass[12pt,a4paper]{article}
\usepackage{amssymb}
\usepackage{amsmath,amsfonts,amsthm}

\usepackage{graphicx}
\usepackage{amsmath}

\usepackage{amsfonts}
\usepackage{amsthm,amscd}
\usepackage{graphicx}
\usepackage{amsmath}
\usepackage{amssymb}
\usepackage{amstext}

\newtheorem{theorem}{Theorem}
\newtheorem{corollary}[theorem]{Corollary}
\newtheorem{definition}[theorem]{Definition}
\newtheorem{example}[theorem]{Example}
\newtheorem{lemma}[theorem]{Lemma}

\newtheorem{proposition}[theorem]{Proposition}
\newtheorem{remark}[theorem]{Remark}

\title{The KMS Condition for the homoclinic equivalence relation and Gibbs probabilities}

\author{A. O. Lopes  and G. Mantovani}

\begin{document}

\maketitle

	
\maketitle

\begin{abstract}

D. Ruelle  considered  a general setting where he is able to characterize equilibrium states  for H\"older potentials based on properties of conjugating homeomorphism in the so called Smale spaces.
On this setting he also shows a relation of KMS states of $C^*$-algebras with equilibrium probabilities of Thermodynamic Formalism. A later paper  by N. Haydn and D. Ruelle  presents a shorter proof of this equivalence.

 Here we consider similar problems but now on the symbolic space  $\Omega = \{1,2,...,d\}^{\mathbb{Z} - \{ 0 \}  }$ and the dynamics will be  given by the shift $\tau$. In the case of  potentials depending on a finite coordinates we will present
a simplified  proof of the equivalence mentioned above which is the main issue of the papers by D. Ruelle and N. Haydn. The class of conjugating homeomorphism
is explicit and reduced to a minimal set of conditions.

We also present  with details (following D. Ruelle) the relation of these probabilities with the KMS dynamical $C^*$-state  on the $C^*$-Algebra associated to the groupoid defined by the homoclinic equivalence relation.

The topics presented here are not new but we believe the main ideas of the proof of the results by Ruelle and Haydn will be quite transparent in our exposition.

\end{abstract}



\section{Introduction} \label{int}

\bigskip

D. Ruelle in  \cite{Reullenoncommutative} considered  a general setting (which includes hyperbolic diffeomorphisms on manifolds)  where he is able to describe a formulation  of the concept of Gibbs state based on {\bf conjugating homeomorphism} in the so called Smale spaces.
On this setting he shows a relation of KMS states of $C^*$-algebras with H\"older equilibrium probabilities of Thermodynamic Formalism. Part of the formulation of this  relation  requires the use of a non trivial  result by N. Haydn (see \cite{HaydnGibbs}). Later, the paper \cite{HD} by N. Haydn and D. Ruelle  presents a shorter proof of the equivalence.

Here we consider similar problems but now on the symbolic space and the dynamics will be  given by the shift. We will present
a simplified  proof of the equivalence mentioned above. The main result of this chapter is Theorem \ref{mai} on section \ref{E}. One can get a characterization of the equilibrium probability for a potential defined on the lattice $\{1,2,...,d\}^{\mathbb{Z} - \{ 0 \}  } $ without using the Ruelle operator (which acts on the lattice $\{1,2,...,d\}^\mathbb{N}  $). The probability we get is invariant for the action of the shift $\tau$ acting on $\{1,2,...,d\}^{\mathbb{Z} - \{ 0 \}  } $.

The proof of this result will take several subsequent sections.

In section \ref{C} we show the relation of these probabilities with the KMS dynamical $C^*$-state  on the $C^*$-Algebra associated to the groupoid defined by the homoclinic equivalence relation. On the initial sections we introduce several results which are necessary for the simplification of the final argument on section \ref{C}.

We present several examples  helping the reader on the understanding of the main concepts.

On \cite{RuelleBook} and also on the beginning of the book \cite{BLL} it is explained the relation of equilibrium states of Thermodynamic Formalism
with the corresponding concept in Statistical Physics. The role of KMS $C^*$-dynamical states on Quantum Statistical Physics is described on \cite{Bra}.  KMS $C^*$-dynamical states correspond to the DLR probabilities (see \cite{CL} for definition) in Statistical Mechanics.

In section \ref{C} we present definitions and properties regarding the $C^*$-algebra we will consider here.

Working on the symbolic space helps to avoid several technicalities which are required in the case of the study of hyperbolic diffeomorphisms on manifolds (where one have to use stable foliation, the local product structure, etc...).

Our proof consider mainly potentials    $A:  \{1,2,...,d\}^{\mathbb{Z} - \{ 0 \}  } \to \mathbb{R}$ which depend on a finite number of coordinates. The case of a general H\"older potential (more technical) can be obtained by adapting our reasoning but we will not address this question here.

On the papers \cite{LM1} and \cite{LO} the authors consider among other things a relation of KMS probabilities with eigenprobabilities for the dual of the Ruelle operator (which are not necessarily invariant for the shift). This problem is analyzed on the lattice $\{1,2,..,d\}^\mathbb{N}$ which is a different setting that the one we consider here. The equivalence relations are also not related.  Despite  some similarities that can be perceived in the statements of the main results obtained in the two settings we point out that the reasoning on the respective proofs are quite different.

Lecture 9 in \cite{GiaII} presents a brief introduction to $C^*$-Algebras and  the KMS condition.

In \cite{EL1} and \cite{EL2} a relation of KMS states in a certain $C^*$-Algebra and eigenprobabilities of the dual of the Ruelle operator is considered.

In a different setting the paper \cite{BK} also considers the homoclinic equivalence relation.

\section{Conjugating homeomorphisms}


In this section $\Omega = \{1,2,...,d\}^{\mathbb{Z} - \{ 0 \}  }$ and a general point $x$ on $\Omega$ is denoted as
$$x= (..., x_{-n},...,x_{-2}  , x_{-1} \,| \,x_1,x_2,...,x_n,...),$$
$x_j \in \{1,2,..,d\}$, $j \in \mathbb{Z}.$

We consider the  dynamics of the shift $\tau: \Omega \to \Omega $, that is,
$$\tau (..., x_{-n},...,x_{-2}  , x_{-1} \,| \,x_1,x_2,...,x_n,...)=(..., x_{-n},...,x_{-2},   x_{-1} , x_1\,|\, x_2,...,x_n,...).$$

We also consider the usual metric $d$ on $\Omega$ which is defined in such way that for $x,y\in \Omega$ we set
$$ d(x,y) = 2^{-N},$$
$N\geq 0$, where for
$$x=
(..., x_{-n}, ... ,x_{-1}\,|\,x_1,..,x_n,..\,)\,\,\,,y=
(..., y_{-n}, ... ,y_{-1}\,|\,y_1,..,y_n,..\,),$$ we have $x_j=y_j$, for all $j$, such that, $-N\leq j\leq N$
and, moreover $x_{N+1}\neq y_{N+1}$, or $x_{-N-1}\neq y_{-N-1}$. Given $x,y$ as above we denote $\vartheta(x,y)=N$, therefore $\vartheta(x,y) = - \log_{2} (d(x,y))$.

Given $x,y \in \Omega$, we say that $x\sim y$ if
$$ \lim_{k \to + \infty} d(\tau^k x, \tau^k y)=0$$
$$\text{and}$$
\begin{equation} \label{equivalencee}
\lim_{k \to - \infty} d(\tau^k x, \tau^k y)=0.
\end{equation}
This means there exists an $N\geq 0$ such that $x_j=y_j$ for $j> N$ and $j< -N$ (note that given $\epsilon>0$, there exists $n$ such that $2^{-n}<\epsilon\leq 2^{-n+1}$, and if $d(x,y)<\epsilon$, then $x$ and $y$ should coincide for coordinates smaller than $n$). In other words, there are only a finite number of $i$'s such that $x_{i} \neq y_{i}$. In this case  we say that $x$ and $y$ are homoclinic.

In this way for large $k>0$ the strings for $\tau^k (x)=z^x$ and $\tau^k (y)=z^y$ are such that $z^x_j= z^y_j$ for $j$
in a large interval $j \in \{-R,-R+1,....,-1,1,....R-1, R\}$, where $R$ is larger with $k$. Then, $ \lim_{k \to + \infty} d(\tau^k x, \tau^k y)=0$.

$\sim$ is an equivalence relation and defines the groupoid $G\subset \Omega  \times \Omega$ of pairs $(x,y)$ of elements which are related (see for instance \cite{Ren0}, \cite{Ren2}, \cite{LM1} or \cite{LO}).

Let $\kappa(x,y)$ be the minimum  $M$ as above. Therefore $x_{\kappa(x,y)} \neq y_{\kappa(x,y)}$ or $x_{-\kappa(x,y)} \neq y_{-\kappa(x,y)}$. Note that $\vartheta(x,y) \leq \kappa(x,y)$ and could be strictly less. Note that $\kappa(x,y)$ is defined just when $x \sim y$.
\medskip

\begin{example}
	For example in $\Omega=\{1,2\}^{\mathbb{Z}-\{0\}}$ take
	$$ x= (...,x_{-n},...,x_{-7},1,2,2,1,2,2\,|\, 1,2,1,2,1,1,x_7,...x_{n},..)$$
	and
	$$ y= (...,y_{-n},...,y_{-7},1,2,2,1,2,2\,|\, 1,2,1,1,1,2,y_7,...y_{n},..)$$
	where $x_j=y_j$ for $|j| > 6=\kappa(x,y).$
	In this case $d(x,y) = 2^{-3}$ and $N=\vartheta(x,y) = 3.$
\end{example}

Given a  H\"older function $U:\Omega \to \mathbb{R}$ it is easy to see that if $x$ and $y$ are homoclinic, then the following function is well
defined
\begin{equation} \label{Vq} V(x,y) = \sum_{n =\,-\,\infty}^\infty (U (\tau^n (x)) -  U (\tau^n (y)) ).\end{equation}

Indded, note that if $x\sim y$, they coincide for large $n$, then, there exists a constant $c$, such that, $d( \tau^n (x), \tau^n (y)) \leq c \,2^{-n}.$ If $U$ has Holder exponent $\alpha$, then, the sum converges absolutely because $\sum_n (2^{\alpha})^{-n} <\infty$.

This function satisfies the property
$$ V(x,y) + V(y,z) = V(x,z)$$
when $x \sim y\sim z$.

A function $V $ with this property will play an important role in some parts of  our reasoning. We will not assume on the first part of this work that $V$ was obtained from a $U$ as above.



\medskip

\vspace{0.5cm}

Now we will describe a certain class of {\bf  conjugating homeomorphism} for the relation $\sim$ (see \eqref{equivalencee}) described above.

Given two fixed points $x$ and $y$ ($y$ in the class of $x$)  we  define the open set $\mathcal{O}_{(x,y)} =B_{\frac{1}{2^{\kappa(x,y)}}}(x)= \{ z \in \Omega : d(x,z) < 2^{-\kappa(x,y)+1}  \} $.

We will define for each such pair $(x,y)$ a conjugating homeomorphisms $\varphi_{(x,y)}$ which has domain on
$\mathcal{O}_{(x,y)}.$

We denote for $m,n \in \mathbb{N}$
$$\overline{ x_{-m}x_{-m+1}...x_{-1}\,|\, x_{1}...x_{n-1}x_{n} }=$$
$$ \{z \in \Omega \,|\, z_j= x_j, \, j=-m, -m+1,...,-1,1,2,...,n-1,n   \},$$

and call it the   cylinder determined by the finite string
$$x_{-m}x_{-m+1}...x_{-1}\,|\, x_{1}...x_{n-1}x_{n}.$$

We will say that a cylinder, or a string, is \textbf{symmetric} if $n=m$.

Note that given   $x\sim y$
$$\mathcal{O}_{(x,y)}=   \overline{ x_{-\kappa(x,y)}\,\, x_{-\kappa(x,y) + 1}\, ... \,x_{-1}\, \vert\,  x_{1}\, ... \,x_{\kappa(x,y) -1}\,\, x_{\kappa(x,y)}},  $$
and $\mathcal{O}_{(x,y)}$ is a symmetric cylinder.

Now we shall define the main kind of {\bf conjugating homeomorphisms} that we will be using. Given $(x,y) \in G$, let $n = \kappa(x,y)$, we define a conjugating $\varphi=\varphi_{(x,y)}$ with domain
$$\mathcal{O}_{(x,y)}=B_{\frac{1}{2^{n}}}(x) = \{ z \in \Omega : d(x,z) < 2^{-n+1}  \}=  \overline{ x_{-n}x_{-n+1}...x_{-1}\,|\, x_{1}...x_{n-1}x_{n} },$$
where $\varphi_{(x,y)}: \mathcal{O}_{(x,y)} \to B_{\frac{1}{2^{n}}}(y)$ is defined by the expression:
$z$ of the form
$$ z= ( ... z_{-n-2} z_{-n-1} \,{\bf \,x_{-n} x_{-n + 1} ... x_{-1} }\vert {\bf x_{1} ... x_{n} \,} \,z_{n+1} z_{n+2} ...)$$
goes to
\begin{equation}\label{varphicorreto}
\varphi_{(x,y)}(z) = ... z_{-n-2} z_{-n-1}\,{\bf  y_{-n} y_{-n + 1} ... y_{-1}
}\vert {\bf y_{1} ... y_{n} }\, z_{n+1} z_{n+2} ...
\end{equation}
We shall call these transformations the family of \textbf{symmetric conjugating homeomorphisms}. We shall denote by $S$ the set of symmetric conjugating homeomorphisms obtained by considering all pairs of related points $x$ and $y$.

\medskip
Note that  the homeomorphism $\varphi_{(x,y)}$ transforms the cylinder $O_{(x,y)} = \overline{ x_{-n}x_{-n+1}...x_{-1}\,|\,x_1 ... x_{n-1}x_{n} }$ in the cylinder $\overline{ y_{-n}y_{-n+1}...y_{-1}\,|\, y_1...y_{n-1}y_{n}}.$

The graph of $\varphi_{(x,y)}$ is on $G$.

A more explicit formulation of the concept of symmetric conjugating homeomorphism  will be presented on next section via expressions  (\ref{kre}) and  (\ref{varphicorreto7}).

\medskip

\medskip
\begin{example}
	Consider
	$$  x= (... 1\,1\,2\,1\,1\,2\,2\,2\,2\,1\,\,\, 1\,1 \vert 2\, 1\, \,\,2\,1\,2\,2\,1\,2\,2\,2\,1\,1...) $$
	and
	$$y =( ... 1\,1\,2\,1\,1\,2\,2\,2\,2\,1\,\,\,1\,2 \vert 1\,2 \, \,\,2\,1\,2\,2\,1\,2\,2\,2\,1\,1 ...) $$
	in this case $\kappa(x,y)=2, $ and for $z$ of the form
	$$ z = (... z_{-4}\,z_{-3}\,\,\, 1\,1 \vert 2\, 1\, \,\,z_3\,z_4\,z_5...)$$
	we get
	
	$$\varphi_{(x,y)}(z) =   (... z_{-4}\,z_{-3}\,\,\, 1\,2 \vert 1\, 2\, \,\,z_3\,z_4\,z_5...).$$

\end{example}

\vspace{0.5cm}

%
%
%
%
%
It is easy to see that the family of symmetric conjugating homeomorphisms we define above has the following properties: given $x \sim y$

a)  $\varphi_{(x,y)} :\mathcal{O}_{(x,y)}\subset \Omega \to \Omega$ is an homeomorphism over its image

b)  $\varphi_{(x,y)}(x)=y$, and

c) $\lim_{k\to \infty} d(\tau^k (z), \tau^k ( \varphi_{(x,y)}(z))=0$ and $\lim_{k\to -\infty} d( \tau^k (z), \tau^k ( \varphi_{(x,y)}(z))=0$.

Item c) implies that $z$ and $\varphi_{(x,y)}(z)$ are on the same homoclinic class.

\medskip


\section{$C^{*}$-Gibbs states and Radon-Nikodym derivative}

We consider the groupoid $G\subset \Omega \times \Omega$ of all pair of points which are related by the homoclinic equivalence relation.

We consider on $G$ the topology generated by sets of the form
$$ \{\,\, (z, \varphi_{(x,y)}(z)) | \,\,\text{where}\,\,z \in \mathcal{O}_{(x,y)}\,\,\text{ with } x \sim y\}.$$

This topology is Hausdorff (see \cite{Reullenoncommutative}).

Now consider a continuous  function $V: G \to \mathbb{R}$ such that
\begin{equation}\label{jairocorrecao2}
V(x,y) + V(y,z) = V(x,z),
\end{equation}
for all related $x,y,z$. Note that this implies that  $V(x,x) = 0$ and $V(x,y) = - V(y,x)$.

Here we call $V$ a modular function.

Under some other notation  the function $\delta(x,y)= e^{ V(x,y)}$  is  called  a modular function (or, a cocycle).

\medskip

\begin{definition} Given a function $V:G \to \mathbb{R}$ as above
	we say that a probability measure $\alpha$ on $\Omega$  is a \textbf{$C^{*}$-Gibbs probability} with respect to the parameter $\beta \in \mathbb{R}$ and $V$, if for any $x\sim y$
	\begin{equation} \label{gibbsR}
	\int_{O_{(x,y)}} \exp(- \beta V(z, \varphi_{(x,y)}(z))) f(\varphi_{(x,y)}(z)) d \alpha(z) =
	\int_{\varphi_{(x,y)}(O_{(x,y)})} f(z) d \alpha (z),
	\end{equation}
	for every continuous function $f : \Omega \to \mathbb{C}$ (and conjugated homeomorphism $(O_{(x,y)}, \varphi_{(x,y)})$).
	
\end{definition}
\medskip

We will show on section \ref{C} a natural relation of this probability $\alpha$ with the $C^*$-dynamical state on a certain $C^*$-algebra. This is the reason for such terminology.

The above definition was taken from \cite{Reullenoncommutative}. This is a version of  the Renault-Radon-Nikodym condition (Def. 1.3.15 in \cite{Ren0}).

\medskip

It is easy to see that the above definition is  equivalent  to say that: given
any pair of finite strings
$$x_{-n} x_{-n+1}...x_{-1}, x_1 ...\,x_{n-1} x_{n}\,\,\text{ and}\,\, y_{-n}y_{-n+1}...y_{-1}\,y_1 ... y_{n-1}y_{n},$$
$n \in \mathbb{N}$, the transformation

\begin{equation} \label{kre} \varphi: \overline{ x_{-n}x_{-n+1}...x_{-1}\,|\, x_1...x_{n-1}x_{n} } \to \overline{ y_{-n}y_{-n+1}...y_{-1}\,|\,y_1...y_{n-1}y_{n}}
\end{equation}
defined by the expression:
\begin{equation}\label{varphicorreto7}
\varphi(z) = (... z_{-n-2} z_{-n-1} \,\,y_{-n} y_{-n + 1} ... y_{-1} \,\vert\, y_{1} ... y_{n}\,\, z_{n+1} z_{n+2} ...),
\end{equation}
where
$$ z=  (... z_{-n-2} z_{-n-1} \,\,z_{-n} z_{-n + 1} ... z_{-1} \,\vert\, z_{1} ... z_{n}\,\, z_{n+1} z_{n+2} ...),$$
is such that for any continuous function $f:  \overline{ y_{-n}y_{-n+1}...y_{n-1}y_{n}} \to \mathbb{R}$
\begin{equation} \label{gibbsR1}
\int_{\overline{ x_{-n}x_{-n+1}...\,|\,...x_{n-1}x_{n} } } e^{- \beta \,V(\,z\,,\, \varphi (z)\,)}  f(\varphi (z)) d \alpha(z) =
\int_{ \overline{ y_{-n}y_{-n+1}...\,|\,...y_{n-1}y_{n}}  } f(z) d \alpha(z).
\end{equation}

\medskip
Note in particulary that by taking $f=1$ we get
\begin{equation} \label{gibbsR2}
\int_{\overline{ x_{-n}x_{-n+1}...\,|\,...x_{n-1}x_{n} } } e^{- \beta \,V(\,z\,,\, \varphi (z)\,)}   d \alpha(z) =
\int_{ \overline{ y_{-n}y_{-n+1}...\,|\,...y_{n-1} y_{n}}}    d \alpha(z).
\end{equation}


In the moment we only  consider symmetric conjugating homeomorphisms of the form $\eqref{varphicorreto7}$.

We will show on section \ref{E} a relation of the \textbf{$C^{*}$-Gibbs probabilities}  $\alpha$ with the {\bf Gibbs (equilibrium) probabilities} of Thermodynamic Formalism.

In a more explicit formulation  $\alpha$ is  such that given any conjugating homeomorphism $(O_{(x,y)}, \varphi_{(x,y)})$ of the form \eqref{kre}, and continuous function $f : \Omega \to \mathbb{C}$
$$
\int_{O_{(x,y)}} e^{-\beta V( z, \varphi_{(x,y)}(z))} f( \varphi_{(x,y)}(z)) d \alpha(z) =$$
$$\int_{O_{(x,y)}} e^{-\beta V( ( ...  z_{-n}, ..., z_{-1} \vert z_{1} ,... ,z_{n},  ...    ), ( ...  z_{-n-1} y_{-n}, ..., y_{-1}  \vert y_{1} ,... ,y_{n} ,z_{n+1},  ...    ))} f(\varphi_{(x,y)}(z)) d \alpha(z) =$$
\begin{equation} \label{jk1} \int_{\varphi_{(x,y)}(O_{(x,y)})} f(z) d \alpha(z).
\end{equation}
In this case, clearly the Radon-Nikodym derivative of the change of coordinates $\varphi $ is
$$  e^{-\beta V( ( ...  z_{-n}, ..., z_{-1} \vert z_{1} ,... ,z_{n},  ...    ), ( ...  z_{-n-1} y_{-n}, ..., y_{-1} \vert y_{1} ,... ,y_{n} ,z_{n+1},  ...    ))}.$$

In order to simplify the notation sometimes on the text we will consider the value $\beta=1$.

\medskip



%
%


\medskip

We will consider a larger class of conjugating homeomorphisms.

\begin{definition}
	Given $n$ and $m$ and
	pair of finite strings
	
	\begin{equation} \label{pre}x_{-n} x_{-n+1}...x_{-1}, x_1 ...\,x_{m-1} x_{m}\,\,\text{ and}\,\, y_{-n}y_{-n+1}...y_{-1}\,y_1 ... y_{n-1}y_{m},
	\end{equation}
	$n,m \in \mathbb{N}$, the transformation
	
	\begin{equation} \varphi: \overline{ x_{-n}x_{-n+1}...x_{m-1}x_{m} } \to \overline{ y_{-n}y_{-n+1}...y_{m-1}y_{m}}
	\end{equation}
	defined by the expression:
	\begin{equation} \label{jairocorrecao}
	\varphi(z) = (... z_{-n-2} z_{-n-1} \,{\bf \,y_{-n} y_{-n + 1} ... y_{-1} \,\vert\, y_{1} ... y_{m}\,}\, z_{m+1} z_{m+2} ...),
	\end{equation}
	where
	$$ z=  (... z_{-n-2} z_{-n-1} \,{\bf \,x_{-n} x_{-n + 1} ... x_{-1} \,\vert\, x_{1} ... x_{m}\,}\, z_{m+1} z_{m+2} ...),$$
	is called a {\bf non-symmetric conjugating homeomorphism} associated to the pair (\ref{pre}).
\end{definition}

Proposition \ref{krt} claims that if $\alpha$ is a $C^{*}$-Gibbs probability, then the relation \eqref{jk1} is satisfied for a bigger class of $\varphi$ transformations, i.e. not necessarily symmetric. Before that  we shall provide the reader with an example of idea of the proof.

\begin{example} \label{wer}
	Consider the non-symmetric conjugating homeomorphism $\varphi : \overline{0 \vert 11} \to \overline{1 \vert 10}$ given by
	$$
	\varphi(... z_{-3} z_{-2} 0 \vert 1 1 z_{3} ...) = ... z_{-3} z_{-2} 1 \vert 1 0 z_{3} ...
	$$
	we shall prove that if $\alpha$ is a $C^{*}$-Gibbs measure then relation \eqref{gibbsR} is valid for $\varphi$. This is actually straightforward, first divide the domain and image of the function into symmetric cylinders,  and in these cylinders apply relation \eqref{jk1}. So in this case consider $\varphi_{0} : \overline{00 \vert 11} \to \overline{01 \vert 10}$, and $\varphi_{1} : \overline{10 \vert 11} \to \overline{11 \vert 10}$ such that
	$$
	\varphi_{a} (... z_{-3} a 0 \vert 1 1 z_{3} ...) =
	(... z_{-3} a 1 \vert 1 0 z_{3} ...)
	$$
	for $a = 0$ or $a = 1$. Now notice that
	$$
	\int_{\overline{0 \vert 11}} e^{-\beta V(x, \varphi(x))} f(\varphi(x)) d \alpha(x) =
	$$
	$$
	\int_{\overline{00 \vert 11}} e^{-\beta V(x, \varphi(x))} f(\varphi(x)) d \alpha(x) + \int_{\overline{10 \vert 11}} e^{-\beta V(x, \varphi(x))} f(\varphi(x)) d \alpha(x) =
	$$
	$$
	\int_{\overline{00 \vert 11}} e^{-\beta V(x, \varphi_{0}(x))} f(\varphi(x)) d \alpha(x) + \int_{\overline{10 \vert 11}} e^{-\beta V(x, \varphi_{1}(x))} f(\varphi(x)) d \alpha(x)
	\stackrel{\eqref{jk1}}{=}
	$$
	$$
	\int_{\overline{01 \vert 10}} f(x) d \alpha(x) + \int_{\overline{11 \vert 10}} f(x) d \alpha(x) =
	\int_{\overline{1 \vert 10}} f(x) d \alpha(x) .
	$$
	This claim proves that relation \eqref{jk1} is valid for this conjugating.
\end{example}

\begin{proposition} \label{krt} Assume $\alpha$ is $C^*$-Gibbs for $V$ as in (\ref{jk1}), then for any
	non-simmetric homeomorphism $(\varphi,\mathcal{O})$, as defined on (\ref{jairocorrecao}), we have that for
	$n,m \in \mathbb{N}$, the transformation
	
	$$
	\int_{ \overline{ x_{-n}x_{-n+1}...x_{-1}\,|\, x_1\,...x_{m-1}x_{m} } } e^{-\beta V( z, \varphi(z))} f( \varphi(z)) d \alpha(z) =$$
	$$\int_{\mathcal{O}} e^{-\beta V( ( ... z_{-n-1} z_{-n}, ..., z_{-1} \vert z_{1} ,... ,z_{m},z_{m+1}  ...    ), ( ...  z_{-n-1} y_{-n}, ..., y_{-1}  \vert y_{1} ,... ,y_{m} ,z_{m+1},  ...    ))} f(\varphi(z)) d \alpha(z) =$$
	\begin{equation} \label{jk3} \int_{\overline{ y_{-n}y_{-n+1}...y_{-1}\,|\, y_1\,...y_{m-1}y_{m}}} f(z) d \alpha(z).
	\end{equation}
	
	\medskip

\end{proposition}

We leave the proof (which is similar to the reasoning of example \ref{wer}) for the reader.

As a particular case we get

\begin{equation} \label{jk09} \int_{ \overline{|\, x_1\,...x_{m} } } e^{-\beta V( z, \varphi(z))} f( \varphi(z)) d \alpha(z) = \int_{\overline{|\, y_1\,...y_{m}}} f(z) d \alpha(z).
\end{equation}
for given  $\overline{|\, x_1\,...\,x_{m} }$, $  \overline{|\, y_1\,...\,y_{m}}$ and the corresponding conjugating homeomorphism $\varphi$.

\medskip
It is possible to consider more general forms of conjugating homeomorphisms as described on the next example.

\begin{example} \label{poi} Consider the homeomorphism $\varphi : \overline{112 \vert 2} \to \overline{1 \vert 122}$ given by
	$$
	\varphi(... z_{-4} 112 \vert 2 z_{2} z_{3} z_{4} ...) =
	(... z_{-4} z_2 z_3 1 \vert 1 2 2 z_4 ...).
	$$
	
	Note that $\overline{112 \vert 2}  $ is translation by $\tau^{-2}$ of the set $\overline{1 \vert 122}.$

	As in the previous example we will prove that if $\alpha$ is a $C^{*}$-Gibbs probability then relation \eqref{jk1} is also valid for such $\varphi$ and $\mathcal{O} = \overline{112 \vert 2}$. First consider the conjugating homeomorphisms, $\varphi_1$, $\varphi_2$, $\varphi_3$ and $\varphi_4$, given by
	$$
	\varphi_1 (... z_{-4} 112 \vert 2\, {\bf 11} z_{4} ...) =
	(... z_{-4}\,{\bf 11} 1 \vert 1 2 2 z_4 ...),
	$$
	$$
	\varphi_2 (... z_{-4} 112 \vert 2\, {\bf 12}\, z_{4} ...) =
	(... z_{-4}\,{\bf  12}\, 1 \vert 1 2 2 z_4 ...),
	$$
	$$
	\varphi_3 (... z_{-4} 112 \vert 2\,{\bf 21}\, z_{4} ...) =
	(... z_{-4} \,{\bf 21}\, 1 \vert 1 2 2 z_4 ...),
	$$
	$$
	\varphi_4 (... z_{-4} 112 \vert 2 \,{\bf 22} z_{4} ...) =
	(... z_{-4} \,{\bf }\,{\bf 22}\, 1 \vert 1 2 2 z_4 ...).
	$$
	
	Therefore we have that
	$$
	\int_{\overline{112 \vert 2}} e^{V(x, \varphi(x))} f( \varphi(x)) d \alpha (x) =
	$$
	$$
	\int_{\overline {112 \vert 211}} e^{V(x, \varphi(x))} f(\varphi(x)) +
	\int_{\overline {112 \vert 212}} e^{V(x, \varphi(x))} f(\varphi(x)) +
	$$
	$$
	\int_{\overline {112 \vert 221}} e^{V(x, \varphi(x))} f(\varphi(x)) +
	\int_{\overline {112 \vert 222}} e^{V(x, \varphi(x))} f(\varphi(x))  =
	$$
	$$
	\int_{\overline {112 \vert 211}} e^{V(x, \varphi_1(x))} f(\varphi_1(x)) +
	\int_{\overline {112 \vert 212}} e^{V(x, \varphi_2(x))} f(\varphi_2(x)) +
	$$
	$$
	\int_{\overline {112 \vert 221}} e^{V(x, \varphi_3(x))} f(\varphi_3(x)) +
	\int_{\overline {112 \vert 222}} e^{V(x, \varphi_4(x))} f(\varphi_4(x)) =
	$$
	$$
	\int_{\overline {111 \vert 122}} f(x) +
	\int_{\overline {121 \vert 122}} f(x) +
	\int_{\overline {211 \vert 122}} f(x) +
	\int_{\overline {221 \vert 122}} f(x) =
	$$
	$$
	\int_{\overline{1 \vert 122}} f d \alpha (x)
	$$
	where some of the $d \alpha$ where omitted. Since we proved that
	$$
	\int_{\overline{112 \vert 2}} e^{V(x, \varphi(x))} f( \varphi(x)) d \alpha (x) =
	\int_{\overline{1 \vert 122}} f d \alpha (x)
	$$
	for any continuous function $f$ then we have that relation (\ref{jk1}) is satisfied.
\end{example}

\vspace{0.5cm}

In analogous way as in last example one can define a conjugating $\varphi$ such that

$\varphi :  \overline{ x_{-n}...{\bf x_{-r}\,...\,x_{-1}}\,|\, x_1\,...x_{m} }\to  \overline{ x_{-n}x_{-n+1}... x_{-r-1}\,|\, {\bf x_{-r}\,... x_{-1}}\,x_1\,...x_{m}}.$

We will consider such transformation $\varphi$  in the next result.
\medskip

\begin{proposition} \label{krtt} Assume $\alpha$ is $C^*$-Gibbs for $V$ as in (\ref{jk1}), then for
	$n,m \in \mathbb{N}$, and $0<r$, such that, $ r\leq n$, we get
	
	$$
	\int_{ \overline{ x_{-n}x_{-n+1}...x_{-r-1}\,{\bf \, x_{-r}\,x_{-r+1}...x_{-1}}\,|\, x_1\,...x_{m-1}x_{m} } } e^{-\beta V( z, \varphi(z))} f( \varphi(z)) d \alpha(z) =$$
	\begin{equation} \label{jk15} \int_{\overline{ x_{-n}x_{-n+1}... x_{-r-1}\,|\, {\bf x_{-r}\,x_{-r+1}... x_{-1}\,x_1\,}\,...x_{m-1}x_{m}}} f(z) d \alpha(z),
	\end{equation}

where $\varphi$ is of the form \eqref{jairocorrecao}.
	
\end{proposition}

{\bf Proof:} The proof is similar to the reasoning of example \ref{poi}. One just has to consider the homeomorphisms
$$\varphi(...z_{-n-r-1}\,  z_{-n-r} ...z_{-n-1}  x_{-n} x_{-n+1}...x_{-1}\,|\, x_1\,...x_{m-1} x_{m} {\bf  z_{m+1}... z_{m+r}}\, z_{m+r+1} ...) =$$
$$
(...  z_{-n-r}\, {\bf z_{m+1}... z_{n+r}\,} x_{-n} x_{-n+1} ... x_{-r-1}\,|\, x_{-r}\,x_{-r+1}......x_{-1} x_{1}...x_{m-1} x_{m}  \, z_{n+r+1}  ...).$$

Note that
$$\tau^{-r} \,(\overline{ x_{-n}x_{-n+1}...x_{-1}\,|\, x_1\,...x_{m-1}x_{m} } \,)=$$
$$ \overline{ x_{-n}x_{-n+1}... x_{-r-1}\,|\, x_{-r}\,x_{-r+1}... x_{-1}\,x_1\,...x_{m-1}x_{m}}.$$
\qed

We want to show that $\alpha$ is $C^*$-Gibbs for $V$, then, the pullback $\rho=\tau^* (\alpha)$ is also $C^*$-Gibbs for $V$.


The next example will help to understand the main reasoning for the proof of the above claim.

\begin{example} Suppose
	$
	V(x,y)
	$ is defined
	when $x \sim y$.  Assume that for all $x,y$ on the groupoid we have that $V(x,y) = V( \tau(x), \tau (y)).$

	Given $\alpha$ consider  the pull back $\rho=\tau^* (\alpha)$.
	
	Consider
	$$\varphi : \overline{11 \vert 2 1}\to \overline{21 \vert 12},$$
	where
	$$ \varphi ( ... x_{-4} x_{-3}\, 1 1 \vert 2 1\, x_3\, x_4...\,) =( ... x_{-4} x_{-3}\, 2 1 \vert 1 2 \, x_3\, x_4...\,),
	$$
	and
	$$\varphi_1 : \overline{112 \vert 1}\to \overline{211 \vert 2},$$
	where
	$$ \varphi_1 ( ... x_{-5} x_{-4}\, 1 1 2 \vert  1\, x_2\, x_3...\,) =( ... x_{-5} x_{-4}\, 2 1 1 \vert  2 \, x_2\, x_3...\,).
	$$
	If for any continuous function $g$ we have that
	$$
	\int_{\overline{11 \vert 21}} e^{V(x, \varphi(x))} g( \varphi(x)) d \alpha (x) =
	\int_{\overline{21 \vert 12}} g (x) d \alpha (x),
	$$
	then, for any continuous function $f$ we have that
	$$
	\int_{\overline{1 1 2 \vert 1}} e^{V(x, \varphi_1(x))} f( \varphi_1(x)) d \rho (x) =
	\int_{\overline{2 1 1 \vert 2}} f (x) d \rho (x).
	$$

	In fact both properties are equivalent.
	
	Note first that $ \varphi_1 \circ \tau = \tau \circ \varphi.$

Moreover, $V(\tau(x), \varphi_1(\tau(x))= V(\tau(x), \tau (\varphi_1(x))= V(x, \varphi_1(x))$ by hypothesis.
	
	Therefore,
	$$
	\int_{\overline{1 1 2 \vert 1}} e^{V(x, \varphi_1(x))} f( \varphi_1(x)) d \rho (x) =$$
	$$\int  I_{\overline{1 1 2 \vert 1}}(x)\,  e^{V(x, \varphi_1(x ))} f( \varphi_1(x)) d \rho (x)=$$
	$$
	\int I_{\overline{1 1 2 \vert 1}}(\tau(x)) e^{V(\tau (x), \varphi_1(\tau(x)))} f( \varphi_1(\tau (x))) d \alpha (x) =
	$$
	$$\int I_{\overline{1 1 2 \vert 1}}(\tau(x)) e^{V(x, \varphi_1(x ))} f( \varphi_1(\tau (x))) d \alpha (x) =
	$$
	$$  \int  I_{\overline{1 1 2 \vert 1}}(\tau(x))\,  e^{V(x, \varphi_1(x ))} f(\tau (\varphi(x)))) d \alpha (x)=
	$$
	$$  \int  I_{\overline{1 1 \vert 2 1}}(x)\,  e^{V(x, \varphi_1(x ))} f(\tau (\varphi(x)))) d \alpha (x)=
	$$
	$$  \int_{\overline{1 1 \vert 2 1}}\,  e^{V(x, \varphi_1(x ))} f(\tau (\varphi(x)))) d \alpha (x)=
	$$
	$$  \int_{\overline{2 1 \vert 1 2}}\,   f(\tau (x)) d \alpha (x)=
	$$
	$$  \int I_{\overline{2 1 \vert 1 2}}\, (x)  f(\tau (x)) d \alpha (x)=
	$$
	$$  \int I_{\overline{2 1 \vert 1 2}}\, (\tau^{-1} \circ \tau ) (x)  f(\tau (x)) d \alpha (x)=
	$$
	$$  \int I_{\overline{2 1 \vert 1 2}}\, (\tau^{-1} (x))  f(x) d \rho (x)=
	$$
	$$  \int_{\overline{2 1 1  \vert 2 }}\,   f(x) d \rho (x).
	$$
	
	Above we took $g = f \circ \tau.$
	
	From the above reasoning we get that  both properties are equivalent.
	
\end{example}

\begin{proposition} \label{kttt}
	If $\alpha$ is $C^*$-Gibbs for $V$, and $V(x,y) = V( \tau(x), \tau (y))$, for all $x,y\in G$, then, the pull back $\rho=\tau^* (\alpha)$ is also $C^*$-Gibbs for $V$.

\end{proposition}

{\bf Proof:} Suppose  $\alpha$ is $C^*$-Gibbs for $V$.

The reasoning of the proof is just a generalization of the argument used on last example.

Consider for $r,s>0$
$$\varphi : \overline{a_{-r}...a_{-1} \vert a_1 a_2...a_s}\to \overline{b_{-r}...b_{-1} \vert b_1 b_2...b_s},$$
where
$$ \varphi ( ... x_{-r+2} x_{-r+1}\, a_{-r}...a_{-1} \vert a_1 a_2...a_s\, x_{s+1}\, x_{s+2}...\,) = $$
$$( ... x_{-r+2} x_{-r+1}\, b_{-r}...b_{-1} \vert b_1 b_2...b_s\, x_{s+1}\, x_{s+2}...\,) ,
$$
and
$$\varphi_1 : \overline{a_{-r}...a_{-1}  a_1 \vert a_2...a_s}\to \overline{b_{-r}...b_{-1}  b_1 \vert b_2...b_s},$$
where
$$ \varphi ( ... x_{-r+2} x_{-r+1}\, a_{-r}...a_{-1}  a_1 \vert a_2...a_s\, x_{s+1}\, x_{s+2}...\,) = $$
$$( ... x_{-r+2} x_{-r+1}\, b_{-r}...b_{-1} b_1 \vert  b_2...b_s\, x_{s+1}\, x_{s+2}...\,) ,
$$

Adapting the argument of last example one can easily show that  if for any continuous function $g$ we have that
\begin{equation} \label{outr}
\int_{\overline{a_{-r}...a_{-1} \vert a_1 a_2...a_s}  } e^{V(x, \varphi(x))} g( \varphi(x)) d \alpha (x) =
\int_{\overline{b_{-r}...b_{-1} \vert b_1 b_2...b_s} } g (x) d \alpha (x),
\end{equation}
then, for any continuous function $f$ we have that
\begin{equation} \label{outr1}
\int_{ \overline{a_{-r}...a_{-1}  a_1 \vert a_2...a_s}} e^{V(x, \varphi_1(x))} f( \varphi_1(x)) d \rho (x) =
\int_{ \overline{b_{-r}...b_{-1}  b_1 \vert b_2...b_s}} f (x) d \rho (x).
\end{equation}

As $\alpha$ is $C^*$-Gibbs for $V$, then (\ref{outr1}) is true for any $f$. From (\ref{outr1}) it follows that
$\rho$ is $C^*$-Gibbs for $V$.

We point out that it is equivalent to ask the $C^*$-Gibbs property for $V$ taking symmetric cylinders or taking not symmetric cylinders (this is implicit on the proof of Proposition \ref{krtt}).

\qed

\bigskip

\section{Modular functions and potentials} \label{mod}

As we mentioned before given a H\"older function $U:\Omega \to \mathbb{R}$ there is a natural way (described by (\ref{Vq})) to get
a continuous function $V$ satisfying the property \eqref{jairocorrecao2}.

We  suppose now that $V$ is such that $ V(x,y) = \sum_{k=-\infty}^\infty [U(\tau^k (x) - U(\tau^k (y)]$, when $x\sim y$, where $U:\Omega \to \mathbb{R}$ is H\"older (see (\ref{Vq})). The function $U$ will sometimes be called a \textbf{potential}. We shall also suppose that $U$ is a finite range potential, or equivalently that it depends on a finite number of positive coordinates, that is,  there is $k\in \mathbb{N}$ and a function $f: \{1,...,d\}^{k} \to \mathbb{R}$, such that, for all $x\in \Omega$ we get
\begin{equation} \label{jairocorrecao3}
U(x)= U \,( ... x_{-n} x_{-n+1}...x_{-2} x_{-1}\,|\, x_1\,x_2...x_{m-1} x_{m} ...)= f(x_1,x_2,...,x_k),
\end{equation}
\smallskip

for this fixed $f$ and $k>0$, where $U:\Omega \to \mathbb{R}$. In this case we say that $U$ depends on $k$ coordinates.

Note that such $V$ satisfies $V(x,y) = V( \tau(x), \tau (y))$ and then Proposition \ref{kttt} can be applied.
\medskip

{\bf Remark 1:} By abuse of language we can write $U: \{1,2,..,d\}^\mathbb{N} \to \mathbb{R}.$
\medskip

If $x \sim y$ it isn't hard to see that there is a finite $M>0$, such that,
$$
V(x,y) =
\sum_{k=-\infty}^\infty [U(\tau^k (x)) - U(\tau^k (y))]=
\sum_{k=-M}^M [U(\tau^k (x)) - U(\tau^k (y))].
$$

In this way, if $z\sim \varphi(z)$, then,
$$ V(z,\varphi(z)) = \sum_{k=-M}^M U(\tau^k (z)) -  \sum_{k=-M}^M U(\tau^k (\varphi(z))= \sum_{k=-M}^M [U(\tau^k (z) - U(\tau^k (\varphi(z))]. $$
Therefore, in this case, equation (\ref{jk1}) means
$$
\int_{\overline{ x_{-n},...,x_{-1}\,|\,x_1,x_2,...,y_n   }} e^{\sum_{k=-M}^M U(\tau^k (\varphi(z))]-\sum_{k=-M}^M U(\tau^k (z)) } f(\varphi(z)) d \alpha(z) =$$
\begin{equation} \label{jk2}\int_{\overline{ y_{-n},...,y_{-1}\,|\,y_1,y_2,...,y_n   }} f(z) d \alpha(z)
\end{equation}
\medskip
If $\alpha$ is $C^{*}$-Gibbs for $V$, and  $ V(z, \varphi(z)) =\sum_{k = - \infty}^{+ \infty} U(z) - U(\varphi(z)) $ we also say {\bf by abuse of language that $\alpha$ is $C^{*}$-Gibbs for $U:\Omega \to \mathbb{R}$}.

\medskip

%
%
%

\vspace{0.5cm}

\begin{definition} Given a function $V:G \to \mathbb{R}$, $
	V(x,y) =
	\sum_{k=-\infty}^\infty [U(\tau^k (x)) - U(\tau^k (y))]$, with $U$ of H\"older class,
	we say that a probability measure $\alpha$ on $\Omega$  is the \textbf{ quasi $C^{*}$-Gibbs probability} with respect to the parameter $\beta \in \mathbb{R}$ and $U$, if there exists constants $d_1>0$ and $d_2>0$, such that, for any $x\sim y$ and any $O_{(x,y)}$,
	$$ d_1 \,\int_{O_{(x,y)}} \exp(- \beta V(z, \varphi_{(x,y)}(z))) g(\varphi_{(x,y)}(z)) d \alpha(z)\,\leq$$
	\begin{equation} \label{qgibbsR}
	\int_{\varphi_{(x,y)}(O_{(x,y)})} g(z) d \alpha (z) \leq  d_2 \,\int_{O_{(x,y)}} \exp(- \beta V(z, \varphi_{(x,y)}(z))) g(\varphi_{(x,y)}(z)) d \alpha(z)
	\end{equation}
	for every every continuous function $g : \Omega \to \mathbb{C}$ (and symmetric conjugated homeomorphism $(O_{(x,y)}, \varphi_{(x,y)}$).
	
\end{definition}
\medskip

In the same way as before one can extend the above property for symmetric conjugated homeomorphisms to non symmetric conjugated homeomorphisms.

\medskip

A $C^*$-Gibbs probability is a quasi $C^*$-Gibbs probability.
\medskip
\medskip

We say that a potential
$\tilde{U}:\{1,2..,d\}^\mathbb{N} \to \mathbb{R}$ - which depends on a finite number of coordinates - is {\bf normalized}, if for $k$ large enough and for any $(x_1,x_2,..,x_k)$ we get $\sum_{j=1}^d e^{\tilde{U}(j,x_1...,x_{k-1})}=1$ - in particularly,
we get $e^{\tilde{U}(x)} =e^{\tilde{U}(x_1...,x_{k})}<1$ for all $x=(x_1,x_2,...)\in \{1,2,...,d\}^\mathbb{N}$.

From this follows that for any $w= (w_1,w_2,...,w_m,...)  \in \{1,2,...,d\}^\mathbb{N}$ and $n\in \mathbb{N}$,
$$\sum_{z_1, z_2,..,z_n =1  }^d e^{ \sum_{j=0}^{n-1} \tilde{U} (\sigma^j(z_1, z_2,..,z_n,w_1,w_2,w_3,...,w_m,...)}=1 $$
where $\sigma$ is the shift acting on $\{1,2,..,d\}^\mathbb{N}.$

Suppose  for such $U$ that $\alpha$ is quasi $C^*$-Gibbs for $U$ (satisfies the double inequality
(\ref{qgibbsR}) for any continuous $g$). This implies in particular that there exist $d_1,d_2>0$, such that, for any cylinders of the form $ \overline{|\,x_1^0,x_2^0,...x_s^0}$ and $\overline{|\,y_1^0,y_2^0,...y_s^0}$,  and a function $\varphi$, such that,
$$d_1
\int_{\overline{|\,x_1^0,x_2^0,...x_s^0} } e^{-\beta V( z, \varphi(z))} g( \varphi_{(x,y)}(z)) d \alpha(z) \leq$$
\begin{equation} \label{jk9}  \int_{\overline{|\,y_1^0,y_2^0,...y_s^0}} g(z) d \alpha(z) \leq d_2
\int_{\overline{|\,x_1^0,x_2^0,...x_s^0} } e^{-\beta V( z, \varphi(z))} g( \varphi_{(x,y)}(z)) d \alpha(z) ,
\end{equation}
where  $\varphi_{(x,y)}$ is the associated conjugating  homeomorphism, such that,
$$\varphi_{(x,y)}: (\overline{|\,x_1^0,x_2^0,...x_s^0} ) \to \overline{|\,y_1^0,y_2^0,...y_s^0} $$

\medskip

\begin{example} \label{poise} Consider the homeomorphism $\varphi : \overline{112 \vert 2} \to \overline{\vert 1122}$ given by
	$$
	\varphi(... z_{-4} 112 \vert 2 z_{2} z_{3} z_{4} ...) =
	(... z_{-4} z_2 z_3 z_4 \vert 1 1 2 2 z_5 ...).
	$$
	
	Note that $\overline{112 \vert 2}  $ is translation by $\tau^{-3}$ of the set $\overline{ \vert 1122}.$

	Consider the conjugating homeomorphisms, $\varphi_1$, $\varphi_2$, $\varphi_3$ and $\varphi_4$, given by
	$$
	\varphi_1 (... z_{-4} 112 \vert 2\, {\bf 11} z_{4} ...) =
	(... z_{-4}\,{\bf 11}  \vert 1 1 2 2 z_5 ...),
	$$
	$$
	\varphi_2 (... z_{-4} 112 \vert 2\, {\bf 12}\, z_{4} ...) =
	(... z_{-4}\,{\bf  12}\,  \vert 1 1 2 2 z_5 ...),
	$$
	$$
	\varphi_3 (... z_{-4} 112 \vert 2\,{\bf 21}\, z_{4} ...) =
	(... z_{-4} \,{\bf 21}\,  \vert 1 1 2 2 z_5 ...),
	$$
	$$
	\varphi_4 (... z_{-4} 112 \vert 2 \,{\bf 22} z_{4} ...) =
	(... z_{-4} \,{\bf }\,{\bf 22}\,  \vert 1 1 2 2 z_5 ...).
	$$
	
	Suppose $\alpha$ is quasi-$C^*$ Gibbs and satisfies (\ref{qgibbsR}).
	
	Therefore,
	$$
	\int_{\overline{112 \vert 2}} e^{V(x, \varphi(x))} f( \varphi(x)) d \alpha (x) =
	$$
	$$
	\int_{\overline {112 \vert 211}} e^{V(x, \varphi(x))} f(\varphi(x)) +
	\int_{\overline {112 \vert 212}} e^{V(x, \varphi(x))} f(\varphi(x)) +
	$$
	$$
	\int_{\overline {112 \vert 221}} e^{V(x, \varphi(x))} f(\varphi(x)) +
	\int_{\overline {112 \vert 222}} e^{V(x, \varphi(x))} f(\varphi(x))  =
	$$
	$$
	\int_{\overline {112 \vert 211}} e^{V(x, \varphi_1(x))} f(\varphi_1(x)) +
	\int_{\overline {112 \vert 212}} e^{V(x, \varphi_2(x))} f(\varphi_2(x)) +
	$$
	$$
	\int_{\overline {112 \vert 221}} e^{V(x, \varphi_3(x))} f(\varphi_3(x)) +
	\int_{\overline {112 \vert 222}} e^{V(x, \varphi_4(x))} f(\varphi_4(x)) \leq
	$$
	$$
	\frac{1}{d_1}\,\,[\,\int_{\overline {11 \vert 1122}} f(x) +
	\int_{\overline {12 \vert 1122}} f(x) +
	\int_{\overline {21 \vert 1122}} f(x) +
	\int_{\overline {22 \vert 1122}} f(x)\,] =
	$$
	$$
	\frac{1}{d_1}\,\int_{\overline{ \vert 1122}} f d \alpha (x),
	$$
	where some of the $d \alpha$ where omitted. We proved that
	$$
	\int_{\overline{112 \vert 2}} e^{V(x, \varphi(x))} f( \varphi(x)) d \alpha (x) \leq \frac{1}{d_1}
	\int_{\overline{\vert 1122}} f d \alpha (x),
	$$
	for any measurable function $f$.
	
	Taking $f=1$, we get that
	$$
	\int_{\overline{112 \vert 2}} e^{V(x, \varphi(x))}  d \alpha (x) \leq \frac{1}{d_1}
	\int_{\overline{\vert 1122}}  d \alpha (x).
	$$

	As $e^{V(x, \varphi(x))}$ is strictly positive we get that if $\alpha(\overline{\vert 1122})=0$, then, $\alpha(\overline{112 \vert 2})=0.$
	
	Using the inequality for $d_2$ in  (\ref{qgibbsR}) we get in a similar way that  if $\alpha(\overline{112 \vert 2})=0$, then, $\alpha(\overline{\vert 1122})=0$.
	
	One can also show that
	$$\int_{\overline{\vert 1122}}  d \alpha (x)\leq d_2
	\int_{\overline{112 \vert 2}} e^{V(x, \varphi(x))}  d \alpha (x) .
	$$
	
\end{example}

\medskip

\begin{proposition} \label{oup} Suppose $\alpha$ is quasi-$C^*$-Gibbs for a \textbf{ potential $U$ that depends on finite coordinates}, then
	$$\alpha( \overline{a_{-r}...a_{-1} \vert a_1 a_2...a_s})>0,$$
	if and only if,
	$$\alpha( \vert \overline{a_{-r}...a_{-1}  a_1 a_2...a_s})>0.$$
	
	Moreover, there exist $b_1>0,b_2>0$, such that, for any  cylinder set of the form  $\overline{a_{-r}...a_{-1} \vert a_1 a_2...a_s}$ we get
	$$   b_1\,\,\,\alpha(\,\overline{a_{-r}...a_{-1} \vert a_1 a_2...a_s}\,) \leq  \alpha(\,\overline{|\,a_{-r}...a_{-1}  a_1 a_2...a_s}\,)\leq$$
	\begin{equation} \label{tutu}
	b_2\,\,\alpha(\,\overline{a_{-r}...a_{-1} \vert a_1 a_2...a_s}\,) .
	\end{equation}

\end{proposition}

{\bf Proof:} We left the proof for the reader which is an adaptation of the reasoning of Example \ref{poise}.


\qed

\medskip
The next result shows that we can always consider normalized potentials (see Theorem 2.2 in \cite{PP} for general results) on the definition of quasi $C^*$-Gibbs  probability.

\begin{theorem} \label{opp} Suppose the probability $\alpha$ on $\Omega$ is $C^{*}$-Gibbs for H\"older potential $U$. Assume, $X:\Omega \to \mathbb{R}$ is such that  $X=U+ g - g \circ \tau+ \lambda$, where $g: \Omega \to \mathbb{R}$ is a H\"older continuous function and $\lambda$ a constant, then $\alpha$ is quasi $C^{*}$-Gibbs for $X$.
\end{theorem}

\textbf{Proof:} Suppose that for any continuous $f$ we have
$$
\int_{O_{x,y}} e^{\beta \sum_{k=-\infty }^\infty U(\tau^k (\varphi(z))) - U(\tau^k (z)) } f(\varphi(z)) d \alpha(z) =$$
\begin{equation} \label{jk7}\int_{\varphi(O_{(x,y)})} f(z) d \alpha(z)
\end{equation}

Note that
$$ \sum_{k=-\infty }^\infty [\, g(\tau^k (z)) -g(\tau^{k} (\varphi(z)))]\,$$
is limited since $g$ is H\"older, actually the summation is absolutely convergent by the same reason. The same can be said of
$$ \sum_{k=-\infty }^\infty [\, g(\tau^{k+1} (z)) -g(\tau^{k+1} (\varphi(z)))]\,$$
\textbf{and of }
$$ \sum_{k=-\infty }^\infty [\,U(\tau^k (z)) - U(\tau^k \varphi(z))\,]$$.

The absolute convergence allow us to sum the quantities above in any order, the resulting sum is limited since each of the above quantities are.

Therefore,
$$ \sum_{k=-\infty }^\infty [ X(\tau^k (\varphi(z)))- X(\tau^k (z))]=$$
$$ [\,\,\sum_{k=-\infty }^\infty U(\tau^k (\varphi(z))]-U(\tau^k (z))\,]+$$
$$ [\,\sum_{k=-\infty }^\infty g(\tau^k (\varphi(z)) - g(\tau^k (z))\,]-$$
$$ [\,\sum_{k=-\infty }^\infty g(\tau^{k+1} (\varphi(z))) -g(\tau^{k+1} (z))\,]$$
is bounded above and below by constants which do not depend on $x\sim y$, $O_{x,y}$ and corresponding $\varphi_{x,y}$.

Then, $\alpha$ is quasi $C^{*}$-Gibbs for $X$.

\qed

By Proposition 1.2 in \cite{PP} given a H\"older potential $U: \Omega \to \mathbb{R}$, one can find $W$ depending on positive coordinates $(1,2,3,..,n,...) \in \{1,2,...,d\}^\mathbb{N}$ and a continuous function $v:\Omega \to \mathbb{R}$ (which depends on finite coordinates),  such that,  $W=U+ v - v \circ \tau.$

The function $V$ is H\"older and then last theorem can be applied.

More precisely, there exist $\tilde{W}:\{1,2,...,d\}^\mathbb{N} \to \mathbb{R}$ an $r$, such that,
$$ W  (...  x_{-n-1} \, \,x_{-n} x_{-n + 1} ... x_{-1} \,\vert\, x_{1} ... x_{m}\, x_{m+1}  ...)=$$
$$ \tilde{W}  ( x_{1} ... x_{m}\, x_{m+1}  ...)=K( x_{1} ... x_{r}\,),$$
for a certain function $K:\{1,2,...,d\}^r \to \mathbb{R}.$

The bottom line is: from Theorem 2.2 in \cite{PP}, given such $\tilde{W}$ one can find, $u$ and positive constant $\lambda$, such that,
$\tilde{W}=\,\tilde{U}+ u - u \circ \tau+ \lambda$. Moreover,  $\tilde{U}: \{1,2,...,d\}^\mathbb{N} \to \mathbb{R}$ and  $u: \{1,2,...,d\}^\mathbb{N} \to \mathbb{R}$ both  depend on a finite number of coordinates.

\medskip

{\bf Remark 2:} Therefore, from Theorem \ref{opp} if $\alpha$ is $C^*$-Gibbs for a H\"older potential $U:\Omega \to \mathbb{R}$, which depends on a finite number of coordinates, we can assume that $\alpha$ is quasi-$C^*$-Gibbs for {\bf another} potential, denoted $\tilde{U}$, which  is normalized and depending on a finite number of coordinates.

By abuse of language one can write  $\tilde{U}: \{1,2,...,d\}^\mathbb{Z} \to \mathbb{R}$.
\medskip

\medskip
\section{Equivalence between equilibrium measures and $C^{*}$-Gibbs measures} \label{E}

First we present two important and well known theorems (see theorems 1.2 and 1.22 in \cite{bowenbook} and also \cite{RuelleBook}).

We will consider without loss of generality that $\beta=1$.

$ \mathcal{M}_{\tau} (\Omega)$ denotes the set on invariant probabilities for $\tau$ acting on $\Omega$.

\begin{theorem} \label{Bowen1.2} (see Theorem 1.2 in \cite{bowenbook})
	Suppose $U : \Omega \to \mathbb{R}$  is of H\"older class.  Then, there is a unique  $\rho \in \mathcal{M}_{\tau} (\Omega)$, for which one can find constants $C_{1} > 0$, $C_{2} > 0$, and P such that, \textbf{for all $s \geq 0$}, for all cylinder $\overline{|\,y_1^0,y_2^0,...y_s^0}$ we have
	\begin{equation} \label{defbowen}
	C_{1} \leq
	\frac{\rho (\overline{|\,y_1^0,y_2^0,...y_s^0})   }
	{\exp \left( -P\, s + \sum_{k=0}^{s-1} U(\tau^{k}x) \right) }
	\leq C_{2},
	\end{equation}
	where
	$$x = (...x_{-k}, x_{-k+1},...,x_{-1}\,|\,x_{1},...,x_m,x_{m+1},...) \in \overline{|\,y_1^0,y_2^0,...y_s^0} \subset \Omega ,$$
	
\end{theorem}

\medskip

We call (\ref{defbowen}) the Bowen's inequalities.

\medskip

\begin{definition} The probability $\rho=\rho_U$ of Theorem \ref{Bowen1.2} is called \textbf{equilibrium probability} for the potential $U$.
\end{definition}

\begin{theorem}
	Given $U$ as above and $\rho_{U}$ the equilibrium measure for $U$, then $\rho_{U}$ is the unique probability  on $\mathcal{M}_{\tau} (\Omega)$, for which
	$$
	h(\rho_U) + \int U d \rho_U = P(U) : =\sup_{\nu \in \mathcal{M}_{\tau}} \{h(\nu) + \int U d \nu\},
	$$
	where $h(\nu)$ is the entropy of $\nu$.
\end{theorem}

For a proof see \cite{PP} or \cite{bowenbook}.

$P(U)$ is called the pressure of  $U$. One can show that the $P$ of (\ref{defbowen}) is equal to such  $P(U).$

\medskip
Remember that if $\alpha$ is $C^{*}$-Gibbs for $V$, and  $ V(z, \varphi(z)) =\sum_{k = - \infty}^{+ \infty} U(z) - U(\varphi(z)) $ we also say by abuse of language that $\alpha$ is $C^{*}$-Gibbs for $U:\Omega \to \mathbb{R}$.

\medskip

Note that if $\rho$ is an equilibrium probability for a H\"older potential $U$, then, it is also an equilibrium probability for
$U + (g \circ \tau)  - g + c$, where $c$ is constant and $g: \Omega \to \mathbb{R}$ is  H\"older continuous (see \cite{PP}). In this way we can assume without lost of generality that $\rho_U$ is an equilibrium probability for a normalized potential $U$. If $U$ is normalized then $P(U)=0$.

If $\alpha$ on $\Omega$ is $C^{*}$-Gibbs for $U$, then, from Remark 2  we have that $\alpha$ is quasi-$C^*$-Gibbs for {\bf another} potential $U$ which is normalized.

Note that given $U$ we are dealing with two  definitions:  $C^{*}$-Gibbs and Equilibrium.  From the above comments  we can assume in either case that $U$ is normalized.

The bottom line is: we can assume (see \cite{PP}) that the H\"older potential $\tilde{U} = U + (g \circ \tau)  - g + c$ is normalized, depends just on future coordinates
$\tilde{U}: \{1,2,...,d\}^\mathbb{N} \to \mathbb{R}$ and has pressure zero.

We will work here (due to Theorem \ref {opp} and the above comments)  with the case where the probability $\alpha$ - which is $C^*$-Gibbs for the potential $U$ - is also a quasi-$C^*$-Gibbs probability for the potential $\tilde{U}$ satisfying Pressure $P(\tilde{U})=0$.  In this case, if we want to prove
expression (\ref{defbowen}) for such probability $\alpha$ over $\Omega$, this can be simplified just showing that there exist $c_1,c_2>0$, such that,
\begin{equation} \label{defbowen23}
c_{1} \leq
\frac{\alpha(\overline{|\,y_1^0,y_2^0,...y_s^0}) }
{\exp \left( \sum_{k=0}^{s-1} \tilde{U}(\sigma^{k}x) \right) }
\leq c_{2},
\end{equation}
where $\sigma$ is the shift acting on $\{1,2,..,d\}^\mathbb{N}$ and
where $x$ is of the form
$$x = (y_1^0,y_2^0,...y_s^0,x_{s+1},...,x_m,x_{m+1},...) \in \{1,2,..,d\}^\mathbb{N} .$$

{\bf Remark 3:} Indeed, due to Remark 2 we get that $\tilde{U} = U + (g \circ \tau)  - g + c$, where $g$ depends on finite coordinates. Therefore,  to show (\ref{defbowen23}) - for {\bf $\alpha$ which is $C^*$-Gibbs for $U:\Omega \to \mathbb{R}$} - is equivalent to prove (see details on the proof of Theorem \ref{opp}) that there exists $C_1,C_2>0$, such that,
\begin{equation} \label{defbowen25}
C_{1} \leq
\frac{\alpha(\overline{|\,y_1^0,y_2^0,...y_s^0}) }
{\exp \left( \sum_{k=0}^{s-1} U(\tau^{k}x) \right) }
\leq C_{2},
\end{equation}
where $\tau$ is the shift acting on $\{1,2,..,d\}^\mathbb{Z}$ and
where
$$x = ( ... x_{-2},x_{-1}\,|\,y_1^0,y_2^0,...y_s^0,x_{s+1},...,x_m,x_{m+1},...) \in \{1,2,..,d\}^\mathbb{Z} .$$

\medskip

It's important to note that the main equivalence (equilibrium and $C^*$-Gibbs)  is still valid in a more general setting of a H\"older potential in a general Smale Space. D. Ruelle proved on the setting of hyperbolic diffeomorphisms  that Equilibrium implies $C^{*}$-Gibbs in his book \cite{RuelleBook}, see theorems 7.17(b), 7.13(b) and section 7.18). On the other hand Haydn proved in the paper \cite{HaydnGibbs} that $C^{*}$-Gibbs implies Equilibrium. Later, the paper \cite{HD} presents a shorter proof of the equivalence.

On the two next sections we will present the proof of the following theorem.

\begin{theorem} \label{mai}
	
	Given a potential $U$ depending on a finite number of coordinates, then,  $\alpha$ is the equilibrium measure for $U$, if and only if, $\alpha$ is
	$C^*$-Gibbs for $U$. As the equilibrium probability is unique we get that the $C^*$-Gibbs probability for $U$ is unique.
	
\end{theorem}

\vspace{7pt}
\section{Equilibrium implies $C^{*}$-Gibbs}


The fact that Equilibrium state  implies $C^{*}$-Gibbs
was proved by Ruelle in a  general setting. The proof is in  the  book  \cite{RuelleBook} (see theorems 7.17(b), 7.13(b) and section 7.18).

For completeness we will explain the proof on our setting.

We drop the $(x,y)$ on
$ \varphi_{(x,y)  }$ and $\mathcal{O}_{(x,y)  }$.

\begin{lemma} \label{7.17Ruelle}
	Let $(\Omega, \tau)$ be the shift on the Bernoulli space $\Omega =\{1,2,...,d\}^{\mathbb{Z}-\{0\}}$ and $\rho_{0}$ be the $\tau$-invariant probability measure which realizes the maximum of the entropy, or, simply the equilibrium state for $U = 0$. If $(\mathcal{O}, \varphi)$ is a conjugating homeomorphism, then for any continuous function $f$
\begin{equation} \label{jairocorrecao4}
	\int_{\mathcal{O}} f(\varphi(x)) d \rho_{0}(x) =
\int_{\mathcal{\varphi(\mathcal{O})}} f(x) d \rho_{0} (x)
\end{equation}

\end{lemma}

{\bf Proof:} Given
$$\mathcal{O}= \overline{ x_{-n}x_{-n+1}...x_{-1}\,|\, x_1\,...x_{m-1}x_{m} } $$ and
$$ \varphi(\mathcal{O})= \overline{ y_{-n}y_{-n+1}...y_{-1}\,|\, y_1\,...y_{m-1}y_{m} },$$
we have that for any $r>m$ and $k>n$
$$\rho_{0} (\overline{ x_{-k}x_{-k+1}...x_{-1}\,|\, x_1\,...x_{r-1}x_{r} }  )  =d^{ -\,(\,r+k\,)}= $$
\begin{equation} \label{vai}  \rho_{0} (\overline{ y_{-k}y_{-k+1}...y_{-1}\,|\, y_1\,...y_{r-1}y_{r} }). \end{equation}

We shall prove that equation \eqref{jairocorrecao4} is valid when $f$ is equal to an characteristic function of an arbitrary cylinder. Note that for this purpose is enough to consider $f$ as the characteristic function of cylinders of the form $\overline{ y_{-k}y_{-k+1}...y_{-1}\,|\, y_1\,...y_{r-1}y_{r} }$.  Therefore,

$$ \int_\mathcal{O} I_{  \overline{ y_{-k}y_{-k+1}...y_{-1}\,|\, y_1\,...y_{r-1}y_{r} }  }(\varphi(x)) d \rho_0(x)=$$
$$     \int_{ \varphi(\mathcal{O}) } I_{  \overline{ y_{-k}y_{-k+1}...y_{-1}\,|\, y_1\,...y_{r-1}y_{r} }  }(y) d \rho_0(y).             $$

From this follows the claim.

The main issue on the above proof is property (\ref{vai}).

\qed

\medskip

We denote by $C^{\alpha}(\Omega)$ the set of $\alpha$ H\"older functions on $\Omega$.

\begin{lemma} \label{7.13Ruelle} (see corollary 7.13 in \cite{RuelleBook})
	Consider the shift space $(\Omega, \tau)$ and $A, B \in C^{\alpha}(\Omega)$. Write for integers $a<0$ and $b>0$
	$$
	Z_{[a,b]} =
	\int\,\,\, e^{ \sum_{k=a}^{b-1} B \circ \tau^{k}   }\,\,d \,\, \rho_{A}
	$$
	Then, $Z^{-1}_{[a,b]} \,\,(\exp \sum_{k=a}^{b-1} B \circ \tau^{k} )\,\, \rho_{A}$ tends to $\rho_{A+B}$ in the weak star topology, when $a \to - \infty$ and $b \to + \infty$.

	In particular, taking $A=0$, when $a \to - \infty$ and $b \to + \infty$, we get that
	$$Z^{-1}_{[a,b]}\,\,\, e^{ \sum_{k=a}^{b-1} B \,\circ\, \tau^{k} }\,\,\, \rho_{0}\,\, \to\,\, \rho_B,$$
	where
	$$
	Z_{[a,b]} =
	\int\,\,\, e^{ \sum_{k=a}^{b-1} B \circ \tau^{k}   }\,\,d \,\, \rho_{0}
	$$
	
\end{lemma}

\begin{theorem} \label{rtw9}
	
	If $\rho_{B}$ is an equilibrium state for a potential $B$ that depends on a finite number of coordinates then it is a $C^{*}$-Gibbs state for $B$.
\end{theorem}

\textbf{Proof:}  The statement holds for $B = 0$ by Lemma \ref{7.17Ruelle}. Moreover, Lemma \ref{7.13Ruelle} allow us to extend this result for all $B \in C^{\alpha} (\Sigma_{N})$ in the following manner: given $\mathcal{O}$ and the associated $\varphi$
$$
\int_{\varphi (\mathcal{O})}  g(x) d \rho_{B} (x) =
\lim_{\substack{a \to - \infty \\ b \to  \infty} } Z^{-1}_{[a,b]} \int_{\varphi(\mathcal{O})} \exp \left( \sum_{k=a}^{b-1} B \circ \tau^{k} (x) \right)  g(x) d \rho_{0} (x) \stackrel{\ref{7.17Ruelle}}{=}
$$
$$
\lim_{\substack{a \to - \infty \\ b \to  \infty}} Z^{-1}_{[a,b]} \int_{\mathcal{O}} \exp \left( \sum_{k=a}^{b-1} B \circ \tau^{k} \circ \varphi (x) \right) g \circ \varphi(x) d \rho_{0} (x) =
$$
$$
\lim_{\substack{a \to - \infty \\ b \to  \infty}} Z^{-1}_{[a,b]}  \int_{\mathcal{O}} \exp \left( \sum_{k=a}^{b-1} B \circ \tau^{k} \circ \varphi (x) - \sum_{k=0}^{b-1} B \circ \tau^{k}(x) \right)
$$
$$
\quad \quad \quad \quad \quad \quad \quad \quad \exp \left( \sum_{k=a}^{b-1} B \circ \tau^{k}(x) \right) g \circ \varphi(x) d \rho_{0} (x) =
$$
$$
\lim_{\substack{a \to - \infty \\ b \to  \infty}} Z^{-1}_{[a,b]} \int_{\mathcal{O}} e^{\left( -V(x, \varphi(x)) \right)} g \circ \varphi(x) \exp \left( \sum_{k=a}^{b-1} B \circ \tau^{k} (x) \right) d \rho_{0} (x) =
$$
$$
\int_{\mathcal{O}} e^{\left( -V(x, \varphi(x)) \right)} g \circ \varphi(x) d \rho_{B} (x).
$$

Since the equality
$$
\int_{\varphi (\mathcal{O})}  g(x) d \rho_{B} (x) =
\int_{\mathcal{O}} e^{\left( -V(x, \varphi(x)) \right)} g \circ \varphi(x) d \rho_{B} (x)
$$
was verified for any conjugating homeomorphism $\varphi$ and any $g$, then it follows that $\rho_{B}$  is an $C^{*}$-Gibbs state for $B$.

\qed

\vspace{7pt}
\section{$C^{*}$-Gibbs implies Equilibrium}

Given a $C^*$-Gibbs probability $\alpha$ for a potential $U$ that depends on a finite number of coordinates we will show in this section that $\alpha$ is the equilibrium probability for $U$. We shall further assume that the potential $U$ depend only on positive coordinates and is normalized according to the Ruelle operator, i.e.
\begin{equation} \label{sav0} \sum_{z_1, z_2,..,z_n =1  }^d e^{ \sum_{j=0}^{n-1} \tilde{U} (\sigma^j(z_1, z_2,..,z_n,w_1,w_2,w_3,...,w_m,...)}=1,\end{equation}
for any $w= (w_1,w_2,...,w_m,...)  \in \{1,2,...,d\}^\mathbb{N}$ and $n\in \mathbb{N}$. Such assuptions aren't restrictive, since given any potential $W$ that depends on a finite number of coordinates, it's possible to find a function $g$ depending on finite coordinates, and a normalized potential $\tilde{W}$ that depends of future coordinates, such that \cite{PP}
$$
W = \tilde{W} + g - g \circ \tau -\lambda
$$





If we show that $\alpha$ is $\tau$-invariant and  also satisfies the Bowen's inequalities for $U$, then, it will follow that $\alpha$ is the equilibrium probability for $U$ by Theorem \ref{Bowen1.2}.

We will show first that a quasi $C^*$-Gibbs probability $\alpha$ for $U$ satisfies the Bowen's inequalities (\ref{defbowen25}) for $U$.

Later we will show that a $C^*$-Gibbs probability $\alpha$ is invariant for $\tau$ (see Proposition \ref{outr2}). This will finally show (see Theorem \ref{rwt3}) that "$C^{*}$-Gibbs implies Equilibrium".

Note that we want to show (\ref{defbowen25}) but due to Remark  3 we just have to show (\ref{defbowen23}).

We assume $\alpha$ is such that (\ref{jk9}) is true, that is, there exists $d_1,d_2>0$, such that, for any continuous function $g$
$$
\,d_1 \int_{\overline{|\,x_1^0,x_2^0,...x_s^0} } e^{- V( z, \varphi_{(x,y)}(z))} g( \varphi_{(x,y)}(z)) d \alpha(z) \leq $$
\begin{equation} \label{jk11}  \int_{\overline{|\,y_1^0,y_2^0,...y_s^0}} g(z) d \alpha(z) \leq d_2 \int_{\overline{|\,x_1^0,x_2^0,...x_s^0} } e^{- V( z, \varphi_{(x,y)}(z))} g( \varphi_{(x,y)}(z)) d \alpha(z).
\end{equation}

\medskip

%
%
%
%

%
%

%


We denote $\mathcal{U}= \sup_{x \in \Omega} U(x) -  \inf_{x \in \Omega} U(x)$.

\begin{lemma} \label{rwt1}

	Given a normalized H\"older potential $U(x) = f(x_1, x_2,..,x_r)$, consider $x_1^0\,...x_s^0$ and  $y_1^0\,...y_s^0$ fixed, and also  $a,b \in \{1,2...,d\}$ fixed. Let
	

	$$x=(...x_{-m}x_{-m+1}...x_{-1}\,|\,  x_1^0\,...x_s^0\, x_{s+1},  x_{s+2},...x_{m-1}x_{m} ...)\in \overline{|\,x_1^0,x_2^0,...x_s^0}  $$
	
	$$y=(...x_{-m}x_{-m+1}...x_{-1}\,|\, y_1^0\,...y_s^0\,x_{s+1}, x_{s+2},...x_{m-1}x_{m} ...) \in \overline{|\,y_1^0,y_2^0,...y_s^0},$$

	and also
	
	$$x_a=(...x_{-m}x_{-m+1}...x_{-1}\,|\,  x_1^0\,...x_s^0\, {\bf a}\, x_{s+2}...x_{m-1}x_{m} ...)\in \overline{|\,x_1^0,x_2^0,...x_s^0}  $$
	
	$$y_b=(...x_{-m}x_{-m+1}...x_{-1}\,|\, y_1^0\,...y_s^0\, {\bf b}\, x_{s+2}...x_{m-1}x_{m} ...) \in \overline{|\,y_1^0,y_2^0,...y_s^0}.$$
	
	Assume that $x \sim y$.

	Then,
	
	$$|\sum_{k = - \infty}^{\infty} U(\tau^{k}(x_a)) - U(\tau^{k} (y_b))|\,\leq $$
	$$2\,r\, \mathcal{U}\,+\,\,|\sum_{k = - \infty}^{\infty} U(\tau^{k}(x)) - U(\tau^{k} (y))|\,.$$
	
\end{lemma}
{\bf Proof:} Let $\mathbb{I}$ the set of indicies for $k$ such that $U(\tau^k (x_a))$ (or, $U(\tau^k (y_b))$) differs from $U(\tau^k (x))$ (or, $U(\tau^k (y))$). It`s easy to see that the cardinality of $\mathbb{I}$ is $r$. Therefore

$$|\sum_{k \in \mathbb{Z}} U(\tau^{k}(x_a)) - U(\tau^{k} (y_b))|\, \leq
$$
$$
|\sum_{k \in \mathbb{Z} \setminus \mathbb{I}} U(\tau^{k}(x_a)) - U(\tau^{k} (y_b))|\,
+ |\sum_{k \in \mathbb{I}} U(\tau^{k}(x_a)) - U(\tau^{k} (y_b))|\ =
$$
$$
|\sum_{k \in \mathbb{Z} \setminus \mathbb{I}} U(\tau^{k}(x)) - U(\tau^{k} (y))|\,
+ |\sum_{k \in \mathbb{I}} U(\tau^{k}(x_a)) - U(\tau^{k} (y_b))|\ \leq
$$
$$
|\sum_{k \in \mathbb{Z} \setminus \mathbb{I}} U(\tau^{k}(x)) - U(\tau^{k} (y))|\ + r \mathcal{U} \leq
$$
$$
|\sum_{k \in \mathbb{Z}} U(\tau^{k}(x)) - U(\tau^{k} (y))|\ + 2r \mathcal{U}
$$

\qed

We will adapt the formulation of Proposition 2.1 in \cite{HaydnGibbs} to the present situation.

\medskip

For fixed
$\overline{|\,x_1^0,x_2^0,...x_s^0} $ denote
$$ U_a = \overline{|\,x_1^0,x_2^0,...,x_s^0, a }$$
$a=1,2..,d.$

Note that  $\sum_a \alpha (U_a) = \,\,\,\alpha (\overline{|\,x_1^0,x_2^0,...x_s^0})$, in particular
\begin{equation} \label{sav1} \sum_a \alpha (U_a) < d \,\,\,\alpha (\overline{|\,x_1^0,x_2^0,...x_s^0}). \end{equation}

Consider now a fixed $\overline{|\,y_1^0,y_2^0,...y_s^0}$ and $\varphi_{a,b}$, $a=1,2,...,d$,  $b=1,2,...,d$, denotes the
conjugating homeomorphism from $U_a$ to  $\overline{|\,y_1^0,y_2^0,...y_s^0\,b}= \varphi_{a,b} (U_a).$

Note also that for each $a$
\begin{equation} \label{sav2}  \alpha (\,\overline{|\,y_1^0,y_2^0,...y_s^0}\,)= \sum_{b=1}^d \alpha (\, \varphi_{a,b}\, (U_a)\,).
\end{equation}

Denote
$$ K = \sup_{m\in \mathbb{N}} \{ \sum_{k=0}^{m-1}    [\tilde{U} \tau^k (u) - \tilde{U} \tau^k (v)]\,\,,\, \text{where}\,\, u,v\in \overline{|\,a_1,a_2,...,a_m}\,\,,$$
$$ \text{and}\,\, (a_1,a_2,...,a_m)\in \{1,2,..,d\}^m\,\}.$$

On the above expression we  ask that $u \sim v$.

Note that if $\alpha$ is $C^*$-Gibbs and satisfies (\ref{jk11}) we get in particular that
$$
\,d_1 \int_{\overline{|\,x_1^0,x_2^0,...x_s^0} } e^{- V( z, \varphi_{(x,y)}(z))}  d \alpha(z) \leq $$
\begin{equation} \label{jk19}  \int_{\overline{|\,y_1^0,y_2^0,...y_s^0}}  d \alpha(z) \leq d_2 \int_{\overline{|\,x_1^0,x_2^0,...x_s^0} } e^{- V( z, \varphi_{(x,y)}(z))} d \alpha(z).
\end{equation}

\medskip
\begin{proposition} \label{este1} Suppose $\alpha$ is quasi-$C^*$-Gibbs for $U$ as above. Then, there exists a constat $c_1>0$, such that,
	
	$$ c_1\leq  e^{ - \sum_{k=0}^{s-1}  U \tau^k (x)}\,\, \alpha(\,\overline{|\,x_1^0,x_2^0,...x_s^0}\,)$$
	for any cylinder $\overline{|\,x_1^0,x_2^0,...x_s^0}$ and any $x$ on the cylinder.

	The $\alpha$-probability of any cylinder is positive.
	
\end{proposition}

{\bf Proof:} We assume that (\ref{jk19}) is true.


Fix a certain cylinder $\overline{|\,x_1^0,x_2^0,...x_s^0}$ and fix a point $x \in \overline{|\,x_1^0,x_2^0,...x_s^0}$ then choose another cylinder $\overline{|\,y_1^0,y_2^0,...y_s^0}$ with non null probability and a point $y \in \overline{|\,y_1^0,y_2^0,...y_s^0}$. Fix $x \in \overline{|\,x_1^0,x_2^0,...x_s^0}$ and $y \in \overline{|\,y_1^0,y_2^0,...y_s^0}$. Choose $a,b \in \{1,2,..., d\}$ and define $x_a$ and $y_b$ as

$$x_a=(...x_{-m}x_{-m+1}...x_{-1}\,|\,  x_1^0\,...x_s^0\, {\bf {a}},  x_{s+2},...x_{m-1}x_{m} ...)\in \overline{|\,x_1^0,x_2^0,...x_s^0}  $$

$$y_b=(...x_{-m}x_{-m+1}...x_{-1}\,|\, y_1^0\,...y_s^0\,{\bf {b}}, x_{s+2},...x_{m-1}x_{m} ...) \in \overline{|\,y_1^0,y_2^0,...y_s^0}.$$

%

we get from Lemma \ref{rwt1} that
$$\alpha (\, \varphi_{a,b}\, (U_a)\,) \leq
d_2 \int_{U_a} e^{\sum_{k=-\infty}^{\infty} U(\tau^k \varphi(z)) - U(\tau^k (z)) } d \alpha (z) \leq
$$
$$
d_2 \int e^{\sum_{k=0}^{s} U(\tau^k \varphi(z)) - U(\tau^k y_b) + U(\tau^k \varphi(z)) -U(\tau^k x_a)} e^{\sum_{k=0}^{s} U(\tau^k y_b) - U(\tau^k x_a) }
$$
$$
e^{\sum_{k=s}^{\infty} U(\tau^k \varphi(z)) - U (\tau^k z)}  e^{\sum_{k = 0}^{\infty} U(\tau^{-k} \varphi(z)) - U (\tau^{-k} z)} d\alpha(z)\leq
$$
$$
d_2 \, e^{ 2 \,K + r\, \mathcal{U} } e^{ \sum_{k=0}^{s-1} [\tilde{U} \tau^k (y_b) - \tilde{U} \tau^k (x_a)]}\,\, \alpha(U_a)\,\leq\,
$$
$$d_2\, e^{ 2\,K + 3\,r\, \mathcal{U}} e^{ \sum_{k=0}^{s-1} [\tilde{U} \tau^k (y) - \tilde{U} \tau^k (x)]}\,\, \alpha(U_a).$$

Then, from (\ref{sav2})

$$ \alpha ( \overline{|\,y_1^0,y_2^0,...y_s^0} )  =  \, \sum_{b=1}^d  \alpha (\, \varphi_{a,b}\, (U_a)\,) \leq d_2\, d\, e^{ 2\,K + 3 r \mathcal{U}} e^{ \sum_{k=0}^{s-1} [U \tau^k (y) - U \tau^k (x)]}\,\, \alpha(U_a) .$$

From this and  from (\ref{sav0}) we get

$$1= \, \sum_{y_1^0,y_2^0,...y_s^0=1 }^d  \alpha ( \overline{|\,y_1^0,y_2^0,...y_s^0} )   \leq d_2\, d\, e^{ 2\,K + 3 \,r \mathcal{U}} e^{ - \sum_{k=0}^{s-1}  U \tau^k (x)}\,\, \alpha(U_a) ,$$
and, finally, for $x =(  ... ,x_{-t},..., x_{-2}, x_{-1}   \,|\, x_1, x_2,..,x_t,...) \in  \overline{|\,x_1^0,x_2^0,...x_s^0}$
$$ d = \, \sum_{a=1}^d \,\sum_{y_1^0,y_2^0,...y_s^0=1 }^d  \alpha ( \overline{|\,y_1^0,y_2^0,...y_s^0} )   \leq \sum_{a=1}^d  d_2\, d\, e^{ 2\,K + 3 r \mathcal{U}} e^{ - \sum_{k=0}^{s-1}  U \tau^k (x)}\,\, \alpha(U_a) =$$
$$ d_2\, d\, e^{ 2\,K + 3\, r \mathcal{U}} e^{ - \sum_{k=0}^{s-1}  U \tau^k (x)}\,\, \alpha(\overline{|\,x_1^0,x_2^0,...x_s^0}).   $$



This also shows that  the $\alpha$-probability of any cylinder $\overline{|\,x_1^0,x_2^0,...x_s^0}$ is positive when $\alpha$ is quasi-$C^*$-Gibbs.

By Proposition \ref{oup} we get that any cylinder of the form $\overline{x_{-m}...x_{-1} \vert x_1 x_2...x_s}$ has positive $\alpha$-probability.

\qed

\medskip

\medskip
\begin{proposition} \label{este2} There exists a constant $c_2>0$, such that,
	
	$$  e^{ - \sum_{k=0}^{s-1}  U \tau^k (x)}\,\, \alpha(\,\overline{|\,x_1^0,x_2^0,...x_s^0}\,) \leq  c_2,$$
	for any cylinder $\overline{|\,x_1^0,x_2^0,...x_s^0}$ and any $x$ on the cylinder.
	
	The $\alpha$-probability of any cylinder is positive.
	
\end{proposition}

{\bf Proof:}  We assume that (\ref{jk19}) is true.


Again consider  fixed $x \in \overline{|\,x_1^0,x_2^0,...x_s^0}$ and $y \in \overline{|\,y_1^0,y_2^0,...y_s^0}$. Choose $a,b \in \{1,2,..., d\}$ and define $x_a$ and $y_b$ as

$$x_a=(...x_{-m}x_{-m+1}...x_{-1}\,|\,  x_1^0\,...x_s^0\, {\bf {a}},  x_{s+2},...x_{m-1}x_{m} ...)\in \overline{|\,x_1^0,x_2^0,...x_s^0}  $$

$$y_b=(...x_{-m}x_{-m+1}...x_{-1}\,|\, y_1^0\,...y_s^0\,{\bf {b}}, x_{s+2},...x_{m-1}x_{m} ...) \in \overline{|\,y_1^0,y_2^0,...y_s^0}.$$

Using an analogous reasoning as in proposition \ref{este2}. But now we use the function $g(z) = e^{V(z, \varphi(z))}$ in the first inequality of \eqref{jk11}. After some algebraic work similar to the former demonstration we reach
$$\alpha(U_a)  \leq \frac{1}{d_1}  \,e^{ 2\,K + r \mathcal{U}} e^{ \sum_{k=0}^{s-1} [U \tau^k (x_a) - U \tau^k (y_b)]}\,\, \alpha( \varphi_{a,b}(U_a)\,)\leq\, $$
$$\frac{1}{d_1}\,e^{ 2\,K + 3\,r\, \mathcal{U}} e^{ \sum_{k=0}^{s-1} [U \tau^k (x) - U \tau^k (y)]}\,\, \alpha (\, \varphi_{a,b}\, (U_a)\,) .$$

Therefore,

$$ e^{ \sum_{k=0}^{s-1}  U \tau^k (y)} \alpha(\overline{|\,x_1^0,x_2^0,...x_s^0})= e^{ \sum_{k=0}^{s-1}  U \tau^k (y)]} \sum_{a=1}^d \alpha(U_a) \leq $$
\begin{equation} \label{aci1} \frac{1}{d_1}\,e^{ 2\,K + 3\,r\, \mathcal{U}} e^{ \sum_{k=0}^{s-1} U \tau^k (x) }\,\, \sum_{a=1}^d \alpha (\, \varphi_{a,b}\, (U_a)\,). \end{equation}

Finally, as $  \overline{|\,y_1^0,y_2^0,...y_s^0\,b}= \varphi_{a,b} (U_a)$ we get  from (\ref{sav0}) and  (\ref{aci1})
$$ d\, \alpha(\overline{|\,x_1^0,x_2^0,...x_s^0})= \sum_{b=1}^d \sum_{y_1^0,y_2^0,...y_s^0=1 }^d  e^{ \sum_{k=0}^{s-1}  U \tau^k (y)} \alpha(\overline{|\,x_1^0,x_2^0,...x_s^0}) \leq $$
$$  \frac{e^{  2 K + 3\, r \mathcal{U}}}{d_1} e^{ \sum_{k=0}^{s-1} U \tau^k (x) }\, \sum_{a=1}^d \,\sum_{b=1}^d \sum_{y_1^0,y_2^0,...y_s^0=1 }^d \alpha (\varphi_{a,b}\, (U_a))= \frac{d\, e^{2 K + 3 \,r \mathcal{U}}}{d_1}\, e^{ \sum_{k=0}^{s-1} U \tau^k (x) }. $$

This shows the claim of the proposition.

\qed

\medskip
Now we have to show that $\alpha$ is invariant by $\tau$.

\begin{corollary} \label{rwt15}
	If $\alpha_1$ and $\alpha_2$ are  quasi $C^*$-Gibbs for $U$,  where
	$$ U (  ... ,x_{-n} ,... ,x_{-2} ,x_{-1} \, | \,x_1,x_2,...,x_r, x_{r+1},...  x_{m} ...   ) = f(x_1,x_2,...,x_r)$$
	for some fixed $r$ and function $f:\{1,2,...,d\}^r\to \mathbb{R}$,
	then $\alpha_1$ is absolutely continuous with respect to $\alpha_2$.
	
\end{corollary}

{\bf Proof:}
We assume that $U$ is normalized. Suppose $\alpha_1$ and $\alpha_2$ are  quasi $C^*$-Gibbs for $U$.

Expression (\ref{defbowen23}) for $\alpha_1$ and $\alpha_2$  will determine, respectively, constants $d_1^1,d_2^1$ and
$d_1^2,d_2^2$.

From last Propositions there exist constants $Y_1>0$ and $Y_2>0$, such that,  for any cylinder $\overline{|\, x_1,x_2,..,x_n}$ and for any point $x$ in this cylinder we get
$$    \frac{\alpha_1 (\overline{|\, x_1,x_2,..,x_n})}{e^{ \sum_{k=0}^{n-1} U(\tau^k (x))}  } \leq Y_1   ,$$
and
$$   Y_2 \leq  \frac{\alpha_2 (\overline{|\, x_1,x_2,..,x_n})}{e^{ \sum_{k=0}^{n-1} U(\tau^k (x))} }.  $$

Therefore,
$$ \frac{Y_2}{Y_1} \alpha_1 (\overline{|\, x_1,x_2,..,x_n})\leq  \alpha_2 (\overline{|\, x_1,x_2,..,x_n}).$$

Now consider a cylinder set of the form
$$(\overline{  x_{-m},...x_{-1}|\, x_1,x_2,..,x_n}).$$

Expression (\ref{tutu}) for $\alpha_1$ and $\alpha_2$ will determine, respectively, constants $b_1^1,b_2^1$ and $b_1^2,b_2^2$.

Then, by  Proposition \ref{oup} we get that

$$b_1^1 \,\alpha_1 (\overline{  x_{-m},...x_{-1}|\, x_1,x_2,..,x_n})\leq
\alpha_1 (\overline{|\,  x_{-m},...x_{-1} x_1,x_2,..,x_n})\leq$$
\begin{equation} \label{uaua} \frac{Y_1}{Y_2} \alpha_2 (\overline{|\,  x_{-m},...x_{-1} x_1,x_2,..,x_n})\leq \frac{Y_1}{Y_2} b^2_2 \alpha_2 (\overline{  x_{-m},...x_{-1}|\, x_1,x_2,..,x_n}).\end{equation}

The Borel sigma-algebra over $\Omega$ is generated by the set of cylinders of the form $\overline{  x_{-m},...x_{-1}|\, x_1,x_2,..,x_n}$.

As the probability $\alpha_j(B)$, $j=1,2$, of a Borel set $B$ is obtained, respectively, as an exterior probability using probabilities of the generators we finally get that the analogous inequalities as in (\ref{uaua}) are true with the same same constants, that is,
\begin{equation} \label{uaua1}    b_1^1 \alpha_1(B) \,\,\,\leq \,\,\,\frac{Y_1}{Y_2} b^2_2 \alpha_2 (B ).\end{equation}

Therefore, $\alpha_1$ is absolutely continuous with respect to $\alpha_2$.

\qed

\begin{proposition} \label{outr2} Assume $\alpha$ is $C^*$-Gibbs for $U$, then,
	$\alpha$ is invariant for $\tau$.

\end{proposition}

{\bf Proof:} From Corollary \ref{rwt15}  we get that any two $C^*$-Gibbs probabilities for $U$ are absolutely continuous with respect to each other.

Suppose $\alpha$ is $C^*$-Gibbs, then, $\alpha_1=\tau^* (\alpha)$ is also $C^*$-Gibbs by Proposition \ref{kttt}. If $\alpha \neq  \tau^* (\alpha)$
then, following Theorem 2.5 in \cite{HD} we get that $\rho_1 = |\alpha_1-\alpha | + \alpha_1-\alpha$ and $\rho_2 = |\alpha_1-\alpha | - \alpha_1 +\alpha$ are also  $C^*$-Gibbs. But $\rho_1$ and $\rho_2$ are singular with respect to each other and this is a contradiction.

Therefore, $\alpha=\tau^* (\alpha)$.

\qed

\begin{theorem} \label{rwt3}

	Suppose $U:\Omega \to \mathbb{R}$ is of the form
	$$ U (  ... ,x_{-n} ,... ,x_{-2} ,x_{-1} \, | \,x_1,x_2,...,x_r, x_{r+1},...  x_{m} ...   ) = f(x_1,x_2,...,x_r),$$
	for some fixed $r$ and fixed function $f:\{1,2,...,d\}^r\to \mathbb{R}$.
	
	If $\alpha$ is  $C^*$-Gibbs for the potential $U$ then $\alpha$ is the equilibrium state for $U$.
	
\end{theorem}

{\bf Proof:} As we know by Proposition \ref{outr2} that $\alpha$ is $\tau$ invariant and, moreover, we also know that $\alpha$ is quasi-$C^*$ invariant for another normalized potential, it follows from  Proposition \ref{este1},  Proposition \ref{este2} and Theorem \ref{Bowen1.2} that $\alpha$ is the equilibrium probability for $U$

\qed

Another conclusion one can get  from the above reasoning is that for potentials that depends on finite coordinates the concepts of quasi $C^*$-Gibbs and $C^*$-Gibbs are equivalent on the lattice $\mathbb{Z}.$

\bigskip

\section{Construction of the $C^{*}$-Algebra} \label{C}

Remember that we consider the groupoid $G\subset \Omega \times \Omega$ of all pair of points which are related by the homoclinic equivalence relation.

Remember also that we consider on $G$ the topology generated by sets of the form
$$ \{\,\, (z, \varphi_{(x,y)}(z))\, | \,\text{where}\,\,z \in \mathcal{O}_{(x,y)}\,\,\text{ and } x,y \in \Omega \text{ such that }x \sim y\}.$$

This topology is Hausdorff \cite{Reullenoncommutative}.

We denote by $[x]$ the class of $x\in \Omega$. For each $x$ the set of elements on the class $[x]$  is countable.


We now come to the construction of the noncommutative algebra. Let $\mathcal{C}_{c} (G)$ be the linear space of complex continuous functions with compact support on $G$. If $A, B \in \mathcal{C}_{c} (G)$ we define the product $A * B$ by
$$
(A*B)(x,y) = \sum_{z \in [x]} A(x,z)B(z,y).
$$
Note that if $(x,y) \in G$ then they are conjugated and so the sum is over all $z$ that are conjugated to $x$ and $y$.

Note that there are only finitely many nonzero terms in the above sum because the functions $A,B$ have compact support \cite{Reullenoncommutative}.

Considering the above, $A*B \in \mathcal{C}_{c}(G)$ as one checks readily, so that $\mathcal{C}_{c}(G)$ becomes an associative complex algebra. An \textbf{involution} $A \to A^{*}$ is defined by
$$
A^{*}(x,y) = \overline{A(y,x)}
$$
where the bar denotes complex conjugation.

For each equivalence class $[x]$ of conjugated points of $\Omega$ there is a representation $\pi_{[x]} \to \mathbb{C}$ in the Hilbert space $l^{2} ([x])$ of square summable functions $[x] \to \mathbb{C}$, such that
$$
((\pi_{[x]}A) \xi (y) = \sum_{z \in [x]} A(y,z) \xi (z)
$$
for $\xi \in l^{2}([x])$. Denoting by $\| \pi_{[x]} A \|$ the operator norm, we write
\begin{equation} \label{supremumnorm}
\| A \| = \sup_{[x]} \| \pi_{[x]} A \|.
\end{equation}

$I_{D}$ (the indicator function of the diagonal $D$) is such that for any $A\in \mathcal{C}_{c} (G)$ we get $I_D * A = A * I_D=A$.

The completion of $\mathcal{C}_{c} (G)$ with respect to this norm is
separable. It is called the    reduced $C^{*}$-algebra which is   denoted by   $C^{*}_r(G)$. The unity element $I_{D}$ is contained in this $C^{*}$ algebra.

\begin{remark} \label{continouscstar}

	If $A \in \mathcal{C}_{c} (G)$ and $t \in \mathbb{R}$, we write
	\begin{equation} \label{invariantstate}
	(\sigma^{t} A)(x,y) = e^{i V(x,y) t} A(x,y)	
	\end{equation}
	defining a one-parameter group $(\sigma^{t})$ of $*$-automorphisms of $\mathcal{C}_{c}(G)$ and a unique extension to a one parameter group of $*$-automorphisms of $C^{*}_r(G)$.
	
	We say that $A \in \mathcal{C}_{c} (G)$ is analytic (a classical terminology on $C^*$-algebras) if the real variable $t$ on the function $t \to \sigma^t A$
	can be extended to the complex variable $z\in \mathbb{C}.$ Under our assumptions this will be always the case. Therefore, $\sigma^{-\beta i} A$ is well defined.
\end{remark}

\begin{definition}
	
	A state $\omega$ on $C^{*}_r(G)$ is a linear functional $\omega: C^{*}_r(G) \to \mathbb{C}$, such that, $\omega(A*A^{*}) \geq 0$, and $\omega (I_D) = 1$ (see \cite{Bra}).
\end{definition}

Such state $\omega$ is sometimes called a dynamical $C^*$-state.

\begin{definition}
	A state $\omega$ is invariant if $\omega \circ \sigma^{t} = \omega$, for all $t\in \mathbb{R}$.
\end{definition}

\medskip

It is of paramount importance to be able to substitute the above real value $t$ by the complex number $\beta \, i$ (where $\beta$ is real). We refer the reader to  Propositions 5.3.6 e 5.3.7 in \cite{Bra} for the technical details of this claim.
\medskip

\begin{definition} \label{jiu} Given a modular function $V: G \to \mathbb{R}$ and the associated  $\sigma_t$, $t \in \mathbb{R},$ we say that
	an invariant state $\omega:C^{*}_r(G)\to \mathbb{C}$ satisfies the \textbf{KMS boundary condition} for $V$ and $\beta\in \mathbb{R}$,  if for all $A, B \in C^{*}_r(G)$, there is a continuous function $F$ on $\{ z \in \mathbb{C} : 0 \leq \text{ Im} (z) \leq \beta \}$, holomorphic in $\{ z \in \mathbb{C} : 0 < \text{ Im} (z) < \beta \}$, and such that for any real $t$
	\begin{equation}\label{KMS}
	\omega(\sigma^{t} A * B) = F(t), \quad \omega(B * \sigma^{t} A) = F(t+i \beta)
	\end{equation}
\end{definition}

\qed

Note that using \eqref{KMS} we have that $F(0) = \omega(A \cdot B)$ and
$$
F(0) =
F(-\beta i + \beta i) =
\omega(B * \sigma^{-\beta i} A).
$$
Therefore, for any $A,B$ we get
$$\omega(A\,*\,B) = \omega(B *  e^{-i \beta} A)$$
which is the classical KMS condition for $\omega$ according to \cite{Bra} (see Propositions 5.3.6 e 5.3.7 there). This condition is equivalent to \textbf{KMS boundary condition}.

\begin{theorem}
	
	If $\mu$ is a probability measure on $\Omega$ then a state $w=\hat{\mu}$ on $C^{*}_r (G)$ can be defined for any $A \in \mathcal{C}_{c} (G)$ by
	
	\begin{equation}
	\hat{\mu}(A) = \int A(x,x) d \mu(x)
	\label{defstate}
	\end{equation}
	
\end{theorem}

{\bf Proof:}

 $\hat{\mu}$ is bounded with respect to the above defined norm.

First note that it's easy to verify  that $\hat{\mu}$ is linear, and for any $A$ we have  $\hat{\mu}(A*A^{*}) \geq 0$ and moreover $\hat{\mu}(I_{D}) = 1$. Now, note that since the diagonal $D$ is a compact set, then any continuous function $A : G \to \mathbb{C}$ has a maximum at $D$, therefore \eqref{defstate} is well defined for continuous function. $\hat{\mu}$ is also well defined on the $C^{*}$-algebra.

\qed

\bigskip

\begin{definition} A probability $\nu $ on $\Omega$  is called a {\bf KMS probability} for the modular function $V$ if the state
	$\hat{\nu}$ on $C^{*}_r (G)$  defined by
	\begin{equation}
	\hat{\nu}(A) = \int A(x,x) \nu(dx)
	\label{defstate1}
	\end{equation}
	satisfies the KMS condition for $V$.  Here $G$ is the groupoid given by the homoclinic equivalence relation.

\end{definition}

This probability is sometimes called quasi-stationary (see \cite{LM1}).

The next claim was proved on \cite{Reullenoncommutative}. For completeness we will present a proof of this claim with full details.

\begin{theorem} If the probability  $\alpha$ on $\Omega$  is a $C^{*}$-Gibbs probability with respect to $V$ and $\beta$, then,
	$\hat{\alpha}$  is a KMS probability for the modular function $\beta \,V$. The associated  $\hat{\alpha}$ is a $C^*$ dynamical state for the $C^*_r (G)$ algebra given by  the groupoid obtained  by the homoclinic equivalence relation and satisfies the KMS boundary condition.

\end{theorem}

{\bf Proof:}
Suppose $\alpha$ is a $C^{*}$-Gibbs state with respect to $\beta V$. We assume $\beta=1$.

$\hat{\alpha}$ is $\sigma^{t}$ invariant if for all $t \in \mathbb{C}$ it's true that
$$
\int \sigma^{t} A(x,x) \alpha (dx) =
\int A(x,x) \alpha (dx)
$$
which by definition \eqref{invariantstate} it's equivalent to
$$
\int e^{iV(x,x)t} A(x,x) \alpha (dx) =
\int A(x,x) \alpha (dx)
$$
but since $V(x,x) = 0$ then the state have to be $\sigma^{t}$ invariant.

Now we will show that if $A, B \in \mathcal{C}_{c}(G)$, then
$$
\hat{\alpha}(\sigma^{t} A * B) =
\int \alpha(dx) \sum_{y \in [x]} e^{iV(x,y)t} A(x,y) \cdot B(y,x)
$$
extends to an entire function (just change $t$ to $z\in \mathbb{C}$). For this purpose we will pick $t_{0} \in \mathbb{C}$ and show that
\begin{equation} \label{basiclimit}
\lim_{t \to t_{0}} \dfrac{\hat{\alpha}(\sigma^{t} A * B) - \hat{\alpha}(\sigma^{t_{o}} A * B)}{t - t_{0}}
\end{equation}
exist. Indeed, the limit \eqref{basiclimit} is equivalent to
$$
\lim_{t \to t_{0}} \dfrac{1}{t - t_{0}} \left( \int \alpha(dx) \sum_{y \in [x]} e^{iV(x,y)t} A(x,y) \cdot B(y,x) - \right.
$$

$$
\left. \int \alpha(ds) \sum_{y \in [s]} e^{iV(s,y)t_{0}} A(s,y) \cdot B(y,s) \right) =
$$
\begin{equation} \label{tempholo}
\lim_{t \to t_{0}} \left( \int \alpha(dx) \sum_{y \in [x]} \dfrac{(e^{iV(x,y)t } - e^{iV(x,y) t_{0}} )}{t - t_{0}} A(x,y) \cdot B(y,x) \right) .
\end{equation}
Always have in mind that for each $x$ the summation is over finite terms.


Let $R$ be a closed ball of radius $1$ centered in $t_{0}$. So we can consider the continuous function $f_{t_{0}} : R \backslash\{t_{0}\} \times $ supp$(A) \to \mathbb{C}$
$$
f_{t_{0}} (t,x) = \sum_{y \in [x]} \dfrac{(e^{iV(x,y)t } - e^{iV(x,y) t_{0}} )}{t - t_{0}} A(x,y) \cdot B(y,x)
$$
To extend $f_{t_{0}}$ for the case $t = t_{0}$ we need to solve the limit
\begin{equation} \label{definitionofL}
L_{t_{0}}(x) =
\lim_{t \to t_{0}} \sum_{y \in [x]} \dfrac{(e^{iV(x,y)t } - e^{iV(x,y) t_{0}} )}{t - t_{0}} A(x,y) \cdot B(y,x) =
\end{equation}
$$
\sum_{y \in [x]} \lim_{t \to t_{0}}  \dfrac{(e^{iV(x,y)t } - e^{iV(x,y) t_{0}} )}{t - t_{0}} A(x,y) \cdot B(y,x) =
$$
$$
\sum_{y \in [x]} iV(x,y) e^{iV(x,y) t_{0}} A(x,y) \cdot B(y,x).
$$
So define $f_{t_{0}} (t_{0},x) = L_{t_0}(x)$.

In this way $f_{t_{0}}$ is a continuous function defined on a compact domain. Therefore we may assume that both it's real and imaginary parts are limited by a value $M$ in the domain. Consider a  sequence of functions indexed by the $t$ variable, $\{f_{t_0}(t_{n}, x)\}_{n \in \mathbb{N^{*}}}$ that converge to $L_{t_0}(x)$ when $n \to \infty$,  e.g. $f_{t_{0}} (t_{0} + (1+i)/n, x)$. In this way the dominated convergence theorem assures that the limit \eqref{basiclimit} is equal to the integral:
$$
\int \alpha(dx) L_{t_{0}}(x).
$$
Indeed formally what we have is,
$$
\int \alpha(dx) L_{t_{0}}(x) =
\int \alpha(dx) \sum_{y \in [x]} iV(x,y) e^{iV(x,y) t_{0}} A(x,y) \cdot B(y,x) =
$$
$$
\int \alpha(dx) \lim_{n \to \infty} \sum_{y \in [x]} \dfrac{(e^{iV(x,y) t_{n}} - e^{iV(x,y) t_{0}} )}{t_{n} - t_{0}} A(x,y) \cdot B(y,x) =
$$
$$
\lim_{n \to \infty} \int \alpha(dx) \sum_{y \in [x]} \dfrac{(e^{iV(x,y) t_{n}} - e^{iV(x,y) t_{0}} )}{t_{n} - t_{0}} A(x,y) \cdot B(y,x) =
$$
\begin{equation} \label{sequencef}
\lim_{n \to \infty} \dfrac{\hat{\alpha}(\sigma^{t_{n}} A * B) - \hat{\alpha}(\sigma^{t_{o}} A * B)}{t_{n} - t_{0}}
\end{equation}
Now since the sequence was arbitrary we could remake these calculations to any desired convergent sequence with the same result, therefore \eqref{sequencef} is equal to
$$
\lim_{t \to t_{0}} \dfrac{\hat{\alpha}(\sigma^{t} A * B) - \hat{\alpha}(\sigma^{t_{o}} A * B)}{t - t_{0}},
$$
what proves existence of the limit in equation \eqref{basiclimit}. This allow us to conclude that $\hat{\alpha}(\sigma^{t} A * B)$ is an holomorphic function everywhere.

Let $F(t) = \hat{\alpha}(\sigma^{t} A * B)$. Using a partition of unity on supp $A$ we may write $A = \sum A_{j}$, where supp $A_{j} \subset W_{j} = \{ (z, \varphi_{j}(z)): z \in \mathcal{O}_{j}\}$, and $(\mathcal{O}_{j}, \varphi_{j})$ is a conjugating homeomorphism.  Since supp $A$ is a compact set then we may assume the summation to occur over a finite amount of elements. Thus

$$
F(t) =\int_{\Omega} \sum_{j} \alpha(dx) A_{j} (x, \varphi_{j} x) B(\varphi_{j} x, x) \exp (i V(x, \varphi_{j} x)t) =
$$
$$
\sum_{j} \int_{\mathcal{O}_{j}} \alpha(dx) A_{j} (x, \varphi_{j} x) B(\varphi_{j} x, x) \exp (i V(x, \varphi_{j} x)t)
$$
and therefore
$$
F(t+ i) =
\sum_{j} \int_{\mathcal{O}_{j}} [ e^{- V(x,\varphi_{j} x)} \alpha (dx) ]
A_{j} (x, \varphi_{j} x) B (\varphi_{j} x, x)  \exp(i V(x, \varphi_{j} x) t)
$$
If $\alpha$ is an $C^{*}$-Gibbs state by \eqref{gibbsR} we have that
$$
F(t+\beta i) =
\sum_{j} \int_{\varphi_{j}(\mathcal{O}_{j})} \alpha (dy) B(y, \varphi^{-1}_{j} y) A_{j} (\varphi^{-1}_{j} y,y) \exp (iV(\varphi^{-1}_{j}y,y)t) =
$$
$$
\sum_{j} \int_{\varphi_{j}(\mathcal{O}_{j})} \alpha (dy) B(y, \varphi^{-1}_{j} y) \sigma^{t} A_{j} (\varphi^{-1}_{j}y,y) =
$$
$$
\int_{\Omega} \sum_{j} \alpha (dy) B(y, \varphi^{-1}_{j} y) \sigma^{t} A_{j} (\varphi^{-1}_{j}y,y) =
$$
$$
\int_{\Omega} \alpha (dy) (B * \sigma^{t} A) (y,y) =
\hat{\alpha} (B*\sigma_{t} A)
$$
so that $\hat{\alpha}$ satisfies the KMS condition.


\qed


\begin{thebibliography}{10}


\bibitem{BLL}
A. Baravieira, R. Leplaideur and A. O. Lopes,
\newblock {\em Ergodic optimization, zero temperature limits and the max-plus
  algebra},
\newblock XXIX Coloq. Bras. Mat  - IMPA - Rio de Janeiro, 2013.


\bibitem{BK}
R. Bissacot and B. Kimura, Gibbs Measures on Multidimensional Sub-
shifts, preprint, 2016.


\bibitem{bowenbook}
R. Bowen,
\newblock {\em Equilibrium States and the Ergodic Theory of Anosov
  Diﬀeomorphisms},
\newblock Springer, 1975.


\bibitem{Bra} O. Bratteli and D. Robinson,
\newblock {\em Operator Algebras and Quantum Statistical Mechanics I,}
\newblock  Springer Verlag.


\bibitem{LM1} G. Castro, A. O. Lopes  and G. Mantovani,
Haar systems, KMS states on von Neumann algebras and $C^*$-algebras on   dynamically defined groupoids and Noncommutative Integration, preprint 2017


\bibitem{CL} L. Cioletti and A. O. Lopes, Interactions, Specifications, DLR probabilities and the Ruelle Operator in the One-Dimensional Lattice, Discrete and Cont. Dyn. Syst. - Series A, Vol 37, Number 12, 6139 -- 6152 (2017)

\bibitem{GiaII} G. DellAntonio,
{\em Lectures on the Mathematics
of Quantum Mechanics II}, Atlantis Press, 2016.


\bibitem{EL1} R. Exel and A. O. Lopes, $C^*$-Algebras, approximately proper equivalence
relations and Thermodynamic Formalism, Erg Theo and Dyn Syst, Vol 24, pp 1051--1082 (2004).


\bibitem{EL2} R. Exel and A. O. Lopes,  C*- Algebras and Thermodynamic Formalism, Sao Paulo Journal of Mathematical Sciences 2, 1 (2008), 285–-307


\bibitem{HaydnGibbs}
N. T. A. Haydn,
\newblock On Gibbs and equilibrium states,
\newblock  Ergod. Th. and Dynam. Sys., (7): 119--132, 1987.


\bibitem{HD}
N. T. A. Haydn and D. Ruelle,
\newblock Equivalence of Gibbs and Equilibrium states for homeomorphisms satisfying expansiveness and specification. \newblock  Comm. in Math. Phys., 148, 155-167, 1992


 \bibitem{KumRen} A. Kumjian and J. Renault, KMS states on   $ C^*$-Algebras associated to expansive maps, Proc. AMS Vol. 134, No. 7, 2067-2078, 2006,




\bibitem{LO} A. O.  Lopes  and E. Oliveira, Continuous groupoids on the symbolic space, quasi-invariant probabilities for Haar systems and the Haar-Ruelle operator,  to appear in  Bull of the Braz. Math. Soc.       (on line 2018)




\bibitem{put1} I. Putnam,
\newblock $C^*$-algebras from Smale spaces.
\newblock  Canad. J. Math. 48, no. 1, 175–-195 , 1996.


\bibitem{PutJ} I. Putnam and J. Spielberg,
\newblock The Structure of $C^*$-Algebras associated with hyperbolic dynamical systems,
\newblock  J. Funct. Anal. 163, no. 2, 279–-299. 1999.


\bibitem{PP}
W. Parry and M. Pollicott,
\newblock {\em Zeta Functions and the Periodic Orbit Structure of Hyperbolic
  Dynamics}.
\newblock 1990.

\bibitem{Pedersen}
G. K. Pedersen,
\newblock {\em $C^{*}$ algebras and Their automorphism groups}.
\newblock Academic Press, 1979.


\bibitem{Ren0}    J. Renault,
\newblock {\em A Groupoid approach to $C^*$-algebras, Lecture Notes in Mathematics 793}
\newblock Springer-Verlag, 1980.





\bibitem{Ren2}
J. Renault,
\newblock{\em $C^*$-Algebras and Dynamical Systems}, XXVII Coloquio Bras. de Matematica
\newblock - IMPA - Rio de Janeiro, 2009




\bibitem{Reullenoncommutative}
D. Ruelle,
\newblock Noncommutative algebras for hyperbolic diffeomorphisms.
\newblock  Inventiones mathematicae, (93):1--13, 1988.


\bibitem{RuelleBook}
D. Ruelle,
\newblock {\em Thermodynamic Formalism}.
\newblock Cambridge University Press, 2 edition, 2004.


\bibitem{Thoms}
C. K. Thomsen,
$C^*$-algebras of homoclinic and heteroclinic structure in expansive dynamics,
Memoirs of Amer. Math. Soc. 206,  No 970, 2010

\end{thebibliography}
\end{document}